\documentclass[12pt]{amsart}
\usepackage{amsthm,amsfonts, amssymb, amscd}
\newcommand{\Q}{{\mathbb Q}}
\newcommand{\Pro}{{\mathbb P}}

\newcommand{\C}{{\mathbb C}}

\newcommand{\m}{{\mathfrak m}}

\DeclareMathOperator{\PGL}{PGL}
\DeclareMathOperator{\Spec}{Spec}
\DeclareMathOperator{\ch}{ch}
\DeclareMathOperator{\Ind}{Ind}
\DeclareMathOperator{\ind}{ind}
\DeclareMathOperator{\res}{res}
\DeclareMathOperator{\GL}{GL}
\DeclareMathOperator{\Td}{Td}
\DeclareMathOperator{\Hom}{Hom}

\newcommand{\liminv}{\lim_{\leftarrow}}

\newtheorem{thm}{Theorem}

\newtheorem{prop}[thm]{Proposition}
\newtheorem{lemma}[thm]{Lemma}
\newtheorem{lemma-definition}[thm]{Lemma-Definition}

\newtheorem{cor}[thm]{Corollary}
\theoremstyle{remark}
\newtheorem{remark}[thm]{Remark}
\newtheorem{example}[thm]{Example}

\newtheorem{definition}[thm]{Definition}
\begin{document}
\numberwithin{thm}{section}
\title[Algebraic cycles and equivariant $K$-theory]
{Algebraic cycles and completions of equivariant $K$-theory}
\author{Dan Edidin}
\address{
Department of Mathematics,
University of Missouri,
Columbia, MO 65211
}
\email{edidin@math.missouri.edu}
\author{William Graham}
\address{
Department of Mathematics,
University of Georgia,
Boyd Graduate Studies Research Center,
Athens, GA 30602
}
\email{wag@math.uga.edu}
\thanks{The second author was supported by N.S.F}

\begin{abstract}
Let $G$ be a complex, linear algebraic group acting
on an algebraic space $X$.
The purpose of this paper is to prove a Riemann-Roch theorem
(Theorem \ref{thm.rrnextgeneration})
which gives a description of the completion of the equivariant 
Grothendieck group $G_0(G,X) \otimes \C$ at any maximal ideal of the
representation ring $R(G) \otimes \C$ in terms of equivariant cycles.
The main new technique for proving this theorem is our non-abelian
completion theorem (Theorem \ref{thm.nonabelian}) for
equivariant $K$-theory. Theorem \ref{thm.nonabelian} generalizes
the classical localization theorems for diagonalizable group actions
to arbitrary groups.
\end{abstract}

\maketitle
\tableofcontents
\section{Introduction}
The now-classic Riemann-Roch theorem of Baum, Fulton and MacPherson
\cite{BFM:75,Ful:84} states that for any separated algebraic space $X$
there is a natural isomorphism between the Grothendieck group $G_0(X)$
of coherent sheaves on $X$ and the Chow group $CH^*X$ of cycles on
$X$.\footnote{Throughout this paper all $K$-groups and
Chow groups are taken with complex coefficients.}  When a linear algebraic group $G$ acts on $X$, the
equivariant Riemann-Roch problem is  to relate the equivariant
Grothendieck group $G_0(G,X)$ of $G$-equivariant coherent sheaves and
the $G$-equivariant Chow group $CH^*_G X$.  In contrast to
the non-equivariant case, elements of $G_0(G,X)$ cannot in
general have unique representations by equivariant algebraic cycles.
For example, if $G$ is a finite group and $X$ is a point then
$G_0(G,X) = \C^r$, where $r$ is the number of conjugacy classes in
$G$, while $CH^*_G(X) = \C$. To obtain precise results we must use the fact
that the equivariant Grothendieck group is a module
for the representation ring $R(G)$.

The main
theorem of \cite{EdGr:00} is an equivariant Riemann-Roch
isomorphism between the completion of $G_0(G,X)$
%\footnote{In \cite{EdGr:00}
%the Grothendieck group of $G$-equivariant coherent sheaves was
%denoted as $G^G(X)$.}
at the ideal in $R(G)$ of virtual representations of rank 0 (the
augmentation ideal) and the infinite product of equivariant Chow
groups $\prod_{n=0}^{\infty} CH^n_G(X)$.  Hence elements of the
augmentation completion of $G_0(G,X)$ may be represented by
equivariant cycles. An obvious problem is to determine whether other
completions of $G_0(G,X)$ admit geometric descriptions.

One of the goals of this paper is to solve this problem for
completions of $G_0(G,X)$ at maximal ideals in $R(G)$ when $G$ is a
complex algebraic group.  Because we are using complex coefficients,
maximal ideals in $R(G)$ correspond bijectively to semisimple
conjugacy classes in $G$.  The correspondence takes a semisimple
conjugacy class $\Psi$ to $\m_{\Psi}$, the maximal ideal of virtual
representations whose characters vanish at $\Psi$;  the augmentation
ideal corresponds to the conjugacy class of the identity element of $G$.
In this paper we generalize the results of \cite{EdGr:00} by
proving that if $h$ is an element of
$\Psi$ with centralizer $Z$, 
then elements of the $\m_\Psi$-adic completion of $G_0(G,X)$
can be represented by $Z$-equivariant cycles on
the fixed subspace $X^h$.  Precisely,
we prove (Theorem \ref{thm.rrnextgeneration}) that for any $G$-space
$X$, there is a natural Riemann-Roch isomorphism $\tau^\Psi_X$ from
the $\m_\Psi$-adic completion of $G_0(G,X)$ to the infinite product
$\prod_{n=0}^\infty CH^n_Z(X^h)$. When $X$ is a smooth scheme, there is
an explicit formula for $\tau^\Psi$ in terms of Chern characters and
Todd classes.

The equivariant Riemann-Roch theorem proved here follows from our
non-abelian completion theorem (Theorem \ref{thm.nonabelian}), a
result of independent interest.  The
progenitor for Theorem \ref{thm.nonabelian} is the classical
localization theorem originally proved by Segal
\cite{Seg:68} for Lie group actions, and extended to algebraic
$K$-theory by Nielsen \cite{Nie:74} and Thomason \cite{Tho:92}.  It
states that if $G$ is a diagonalizable group (e.g., a torus) and $h$
is an element of $G$, then the pushforward $G(G,X^h) \stackrel{i_*}
\to G(G,X)$ becomes an isomorphism after localizing at the maximal
ideal $\m_h$ of $R(G)$. (Here $G(G,X)$ denotes the infinite direct sum
of equivariant $K$-groups $\oplus_{n = 0}^\infty G_n(G,X)$.)   
Moreover, if $X$ is smooth then there is an explicit
localization formula
\begin{equation} \label{eqn.classical} 
\alpha = i_* \left( {i^*\alpha\over{\lambda_{-1}(N_i^*)}}\right)
\end{equation}
where $\alpha \in G(G,X)_{\m_h}$ and $N_i^*$ is the conormal
bundle of the regular embedding $i \colon X^h \to X$. 
Many applications of the localization
theorem, such as a simple proof of the Weyl character formula, come 
from this explicit localization formula.
In \cite{EdGr:05}, we extended this
result to the case where $G$ is an arbitrary algebraic group, but with
the assumption that $G$ acts with finite stabilizer, and used our version
of the explicit localization formula 
to give a Riemann-Roch formula for orbifolds. Related results
were also obtained in that case by Toen \cite{Toe:99} and Vezossi-Vistoli
\cite{VeVi:02}.

In this paper, we remove the restriction on the group action, but to
do so, we must complete rather than localize.  (In the finite
stabilizer case these two operations coincide because the equivariant
$K$-theory is supported at a finite number of maximal ideals in
$R(G)$.)  The nonabelian completion theorem states that there is a
pushforward isomorphism
$$
i_!: \widehat{G(Z,X^h)}_h \to \widehat{G(G,X)},
$$
where $\widehat{G(Z,X^h)}_h$ is the completion of $G(Z,X^h)$ at the
maximal ideal $\m_h$ of $R(Z)$, and $\widehat{G(G,X)}$ is the
completion of $G(G,X)$ at $\m_{\Psi}$.  Moreover, when $X$ is smooth
a formula analogous to \eqref{eqn.classical}
above holds. This formula allows us to obtain our explicit formula
for the Riemann-Roch isomorphism $\tau^\Psi_X$ when $X$ is smooth.

The proof of the nonabelian completion theorem rests on the
construction of a ``twisted induction" functor in equivariant
$K$-theory.  If $G$ is a connected reductive group and $Z \subset G$
is a connected subgroup of maximal rank such that both
groups have simply connected commutator subgroups then we can use a theorem of
Merkurjev to define an induction map $\ind \colon G(Z,X) \to G(G,X)$
satisfying several natural properties, including what we call a
``twisted reciprocity" formula relating $\ind$ to the usual restriction functor
$\res \colon G(G,X) \to G(Z,X)$.  When $Z$ is the centralizer of a
semi-simple element $h \in G$, generalities about completions imply
that the $\m_h$-completion $\widehat{G(Z,X)}_h$ of $G(Z,X)$ is
naturally identified as a summand in the $\m_{\Psi}$-adic completion
$\widehat{G(Z,X)}$.  By composing the inclusion as a summand with the
map $\ind$, we obtain a map
$$
\ind_h \colon \widehat{G(Z,X)}_h \to \widehat{G(G,X)}.
$$  Using some basic facts about invariants and completion,
we prove that $\ind_h$ and the natural map 
$$
\res_h \colon \widehat{G(G,X)} \to  \widehat{G(Z,X)}_h
$$ are inverse isomorphisms (Proposition \ref{prop.indres}).  (This
proposition illustrates the necessity of working with completions
rather than localizations, as the natural restriction map of localizations
$G(G,X)_{\m_\Psi} \to G(Z,X)_{\m_h}$ is not in general an isomorphism
(Example \ref{ex.needcompletions}).)  In the case that $Z$ and $G$ are
both connected and reductive and $G$ has simply
connected commutator subgroup, the map $i_!$ can be defined as a
composition
$$
\widehat{G(Z,X^h)}_h \stackrel{i_*} \to \widehat{G(Z,X)}_h \stackrel{\ind_h}
\to \widehat{G(G,X)}.
$$ Because $h$ is in the center of $Z$, we can argue as in the proof
of the localization theorem for tori to show that $i_*$ is an
isomorphism; hence $i_!$ is an isomorphism as well.  For general $G$
and $Z$ a change of groups argument similar to that of \cite{EdGr:05}
can be used to define $i_!$.

There are a number of natural questions arising from this work.   
If the element $h$ is defined over a subfield
$L \subset \C$ then our techniques imply that Theorem \ref{thm.nonabelian}
holds for completions of equivariant $K$-theory tensored with $L$. 
An interesting problem for further study is to prove a version of
Theorem \ref{thm.nonabelian} for completions of integral or rational
equivariant $K$-theory at maximal (or prime) ideals of the representation
ring $R(G)$. It would also be interesting to extend these results to
algebraic groups defined over arbitrary fields.  Similarly, if $\Psi$ is the conjugacy
class of an element $h \in G$ which is defined over a subfield $L \subset \C$
then our techniques show that there is an isomorphism of completions
of equivariant $K$-theory tensored with $L$ and equivariant Chow groups with coefficients
in $L$. Given an arbitrary maximal (or prime) ideal ${\mathfrak m}$ in $R(G) \otimes \Q$
a very interesting open problem is to represent elements of the ${\mathfrak m}$-adic
completion of $G_0(G,X)$ as formal equivariant algebraic cycles on an appropriate
subspace of $X$.
Another natural question is whether a version of Theorem \ref{thm.rrnextgeneration} holds
for higher $K$-theory. Specifically one can ask if there is a natural
isomorphism between the ${\mathfrak m}_\Psi$-adic completion of $G_n(G,X)$
and the infinite product of higher equivariant Chow groups $\prod_{i=0}^\infty CH^i_G(X,n)$.
Such a result would follow from the methods of this paper and the higher
$K$-theory version of the Riemann-Roch isomorphism of \cite{EdGr:00}. Unfortunately
we have not been able to construct such an isomorphism due to the difficulties
in comparing different completions of higher equivariant $K$-theory (cf \cite[Remark 2.2]{EdGr:00}).

{\bf Acknowledgement:} The authors are very grateful
to the referees for their careful reading as well as for many comments
which helped to clarify the exposition.

\subsection{Background} \label{ss.background}
We work entirely over the ground field $\C$ of complex numbers. All
algebraic spaces are assumed to be separated and of finite type over
$\C$. For a reference on the theory of algebraic spaces see the book
\cite{Knu:71}.  As in \cite{EdGr:98} we will refer to an integral
subspace of an algebraic space $X$ as a {\it subvariety} of $X$.  Note
that if $h$ is an automorphism of a separated algebraic space $X$,
then the fixed point subspace $X^h$ is closed in $X$, since it is the
inverse image of the diagonal under the morphism $X \to X \times X$
given by the graph of $h$.

All algebraic groups are assumed to be linear. A basic reference for
the theory of algebraic groups is \cite{Bor:91}.  If $G$ is an algebraic
group and $h \in G$, then ${\mathcal Z}_G(h)$ denotes the centralizer
of $h$ in $G$, and $C_G(h)$ the conjugacy class of $h$ in $G$.

Throughout the paper we take the complex numbers as our coefficients 
for $K$-groups and Chow groups, e.g., we write $G_0(G,X)$ for
$G_0(G,X) \otimes \C$.

\subsubsection{Representation rings}
We recall here some facts about representation rings proved in \cite{EdGr:05}.
If $G$ is an algebraic group then $R(G)$
denotes the representation ring of $G$ tensored with $\C$.
By \cite[Proposition 2.5]{EdGr:05} there is a bijective
correspondence between semisimple conjugacy classes in $G$
and maximal ideals in $R(G)$. If $\Psi$ is the conjugacy 
class of a semisimple element, then the corresponding maximal
ideal $\m_\Psi$ consists of virtual representations whose
characters vanish at $\Psi$. 

If $G$ is an algebraic group and $H \subset G$ is a closed subgroup
then the restriction map $R(G) \to R(H)$ is a finite morphism
(\cite[Proposition 2.3]{EdGr:05}).  As a result, if $\Psi$ is a semisimple conjugacy
class in $G$, then $\Psi \cap H$ decomposes into a finite number of
conjugacy classes $\Psi_1, \ldots \Psi_l$, and a maximal ideal
$\m_{\Psi'} \subset R(H)$ lies over $\m_\Psi \subset R(G)$ if and only
if $\Psi' = \Psi_k$ for some $k$ (\cite[Proposition 2.6]{EdGr:05}).

\subsubsection{Equivariant $K$-theory}
As in \cite{EdGr:05}, $G(G,X)$ denotes the infinite
direct sum $\oplus_{i=0}^\infty G_i(G,X),$ where $G_i(G,X)$ is the
$i$-th Quillen $K$-group of the category of $G$-equivariant coherent
sheaves on $X$, tensored with $\C$. Thus, $G_0(G,X)$ is the
Grothendieck group of $G$-equivariant coherent sheaves, tensored with
$\C$. Likewise, $K_0(G,X)$ denotes the
Grothendieck ring of $G$-equivariant vector bundles, also tensored
with $\C$.
When $X$ is a smooth scheme then $K_0(G,X)$ and $G_0(G,X)$
may be identified thanks to Thomason's equivariant resolution theorem
\cite{Tho:87}. If $X = \Spec \C$ then $K_0(G,X) = G_0(G,X) = R(G)$.

If $N$ is a $G$-equivariant vector bundle of rank $n$ on $X$ then
$\lambda_{-1}(N)$ is the formal sum $\sum_{i=0}^{n} (-1)^i[\Lambda^{i}N]$.

If $Z$ is a closed subgroup of $G$ and $X$ is a $Z$-space, then
there is a ``Morita equivalence" identification of $G(Z,X)$ with $G(G,G \times_Z X)$
(see Section 3.1 of \cite{EdGr:05} for more details).  In particular,
an element of $R(Z)$ such as $\lambda_{-1}({\mathfrak g^*}/{\mathfrak z^*})$
can be viewed as operating on $G(G,G \times_Z X)$.
If $X$ is a $G$-space, then $G \times_Z X$ is isomorphic to $G/Z \times X$,
and the relative tangent bundle of the morphism $G \times_Z X \to X$
is identified under this equivalence with the pullback
of the element ${\mathfrak g}/{\mathfrak z}$ of $R(Z)$ to $G(Z,X)$.  In this case
we will often simply identify the relative tangent bundle of the morphism as
${\mathfrak g}/{\mathfrak z}$ without further comment.  

For a sketch of some basic properties of equivariant
$K$-theory, most of which are due to Thomason, we refer the reader to
Section 3.1 of our paper \cite{EdGr:05}. Facts from that section will
be used in this paper without further reference.

\section{Induced representations in algebraic $K$-theory}
Let $G$ be a connected reductive group and let $H \subset G$ be a
connected reductive subgroup of the same rank. Assume that $G$ 
and $H$ have simply connected commutator subgroups.
The goal of this
section is to define a natural twisted induction map $G(H,X) \to
G(G,X)$ satisfying an appropriate reciprocity formula
(\ref{prop.weylinduction}). A theorem of Merkurjev is crucial
for proving properties of the twisted induction map.

\subsection{Merkurjev's theorem and Weyl group actions} \label{ss.subgroups}
Suppose that $H$ is a subgroup of $G$ and $X$ is a $G$-space.  There
is a natural map
\begin{equation} \label{e.merkurjev}
\theta: R(H) \otimes_{R(G)} G(G,X) \rightarrow G(H,X),
\end{equation}
given by $\theta([V] \otimes \alpha) = [V] \res \alpha$.  This
map can be realized another way (cf. \cite[Prop.~3.2]{EdGr:00}).
We can identify $G \times_H X$
with $G/H \times X$ via the map taking $(g,x) H$ to $(gH,gx)$.
This identification is $G$-equivariant, where $G$ acts on $G \times_H
X$ by left multiplication on $G$, and on $G/H \times X$ by the direct
product of the actions on $G/H$ and $X$.  Consequently
the category ${\tt coh}_X^H$ of $H$-equivariant coherent sheaves
on $X$ is equivalent to the category ${\tt coh}_{G/H \times X}^G$
of $G$-equivariant coherent sheaves on $G/H \times X$.
In terms of this equivalence the map $\theta$ is given as follows.
Let $\pi_1$ and $\pi_2$ denote the projections of $G/H \times X$ to
$G/H$ and $X$, respectively, and given a representation $V$ of $H$,
let ${\mathcal V}$ denote the vector bundle $G \times_H V$ on $G/H$.
Then $\theta$ is the map $R(H) \otimes_{R(G)} G(G,X) \to 
G(H,X)$ induced by the family of functors 
$f_V \colon {\tt coh}_G X \to {\tt coh}_H X$,
${\mathcal F} \mapsto \pi_1^*{\mathcal V} \otimes \pi^*_2 {\mathcal F}$,
defined for each $H$-module $V$.
We may write $\theta_H^G$ for $\theta$ if we wish to make explicit the
groups involved.

If $G$ is a connected reductive group with simply connected commutator subgroup and $T$ is a
maximal torus of $G$, then Merkurjev proves (\cite[Proposition 31]{Mer:05})
that the map
$\theta_T^G \colon R(T) \otimes_{R(G)} G(G,X) \to G(T,X)$ is an isomorphism.
(Note that Merkurjev's result holds for $K$-theory with integral coefficients.)
The following 
extension of Merkurjev's result will be useful for proving properties of the
twisted induction map.
\begin{prop} \label{prop.tensor}
Suppose that $H \subset G$ are connected reductive groups of the same
rank both with simply connected commutator subgroups. 
Then the natural map $\theta_H^G \colon R(H) \otimes_{R(G)} G(G,X) \to G(H,X)$
taking
$[V] \otimes \alpha$ to $\res(\alpha)[V]$ is an isomorphism.
\end{prop}

\begin{proof}
Let $T$ be a maximal torus of $H$ (hence also of $G$).  By Merkurjev's theorem,
$$G(T,X) =   R(T) \otimes_{R(G)} G(G,X)= R(T) \otimes_{R(H)}G(H,X) .$$
Hence the map
$$
R(T) \otimes_{R(H)} \left(R(H) \otimes_{R(G)} G(G,X)\right) \to
R(T) \otimes_{R(H)} G(H,X)
$$
is an isomorphism, since both sides are
identified with $G(T,X)$.  Since $H$ is connected and reductive, the
restriction map $R(H) \to R(T)$ is a split injection. Thus tensoring with $R(T)$ is a fully faithful
functor. Therefore $G(G,X) \otimes_{R(G)} R(H) =
G(H,X)$ as claimed.
\end{proof}

Let $G$ be as above; and let $T$ be a maximal torus
of $G$ and $W$ the corresponding Weyl group.  If $X$ is a
$G$-space, then $W$ acts 
on $G/T \times X$ by the rule $w \cdot (gT, x) = (gw^{-1}T, x)$.
This action commutes with the action of $G$ and hence induces
an action of $W$ on $G(G, G/T \times X)$, which we can
identify with $G(T,X)$.  This action of $W$ is natural with
respect to flat pullbacks and proper pushforwards arising
from morphisms of $G$-spaces.   The isomorphism
$$
\theta_T^G: R(T) \otimes_{R(G)} G(G,X) \rightarrow G(T,X)
$$
is $W$-equivariant, where $W$ acts on the source by its action
on $R(T)$.  By a theorem of Serre \cite[Theorem 4]{Ser:68} we may identify $R(G) = R(T)^W$, so $G(T,X)^W = G(G,X)$.
(Note that Serre's result holds even if the commutator subgroup is not simply
connected.)

\subsection{Twisted induction}
\begin{definition}
Let $H \subseteq G$ be connected reductive groups of equal rank with simply 
connected commutator subgroups, and
let $T$ be a maximal torus of $H$ (hence also of $G$).  Let $W_1
\subseteq W$ denote the Weyl groups of $T$ in $H$ and $G$,
respectively.  Let $X$ be a $G$-space.  We define
$$
\ind^G_H: G(H,X) \to G(G,X)
$$ as follows.  Identify $G(G,X)$ (resp.~$G(H,X)$) as the
$W$-invariants (resp.~$W_1$-invariants) in $G(T,X)$, and for any
$\alpha \in G(H,X)$, set
$$
\ind^G_H(\alpha) = \sum_{w \in W/W_1} w \alpha.
$$
(Note that we use the same notation for elements of $W/W_1$ and
lifts to $W$; sometimes we also use the same notation for
lifts of elements of $W$ to $G$.)
\end{definition}

\begin{definition}
If $G$ is a reductive group, define a linear map
$R(G) \to \C$, $\alpha \mapsto \alpha^G$, by setting
$[V]^G = \dim V^G$ for any $G$-module $V$ and extending 
by linearity. Likewise, define a symmetric bilinear map
$\Hom_G(\_, \_)\colon R(G) \otimes R(G) \to \C$ by setting
$$\Hom_G([V], [W]) = [V^* \otimes W]^G$$ for $G$-modules
$V, W$.
\end{definition}

The twisted induction map has properties given in the following proposition.  
Part (d) is analogous to the Frobenius reciprocity formula for induced
representations, and justifies our use of the term twisted induction map.

\begin{prop} \label{prop.weylinduction}
Let $H \subseteq G$ be connected reductive groups of equal rank.  Assume
that both $H$ and $G$ have simply connected commutator subgroups. The map
$$
\ind^G_H: G(H,X) \to G(G,X)
$$
has the following properties.

(a) $\ind^G_H$ is natural with respect to flat pullbacks and proper pushforwards
arising from morphisms of $G$-spaces.

(b) If $K \subset H \subset G$ are connected reductive groups of equal rank
with simply connected commutator subgroups, then
$$
\ind^G_H \circ \ind^H_K  = \ind^G_K.
$$

(c) (Projection formula) If $\alpha \in R(H)$ and $\beta \in G(G,X)$ then
$$
\ind^G_H(\alpha \res (\beta) ) = \ind^G_H(\alpha) \beta.$$

(d) (Reciprocity formula) When $X = \Spec \C$, the map $\ind^G_{H} \colon R(H) \to R(G)$
satisfies
\begin{equation} \label{eq.recip}
\Hom_G(\ind^G_{H}\alpha, \beta) = \Hom_H( \lambda_{-1}({\mathfrak g^*}/
{\mathfrak h^*})\alpha, \res \beta),
\end{equation}
where $\mathfrak g$ and $\mathfrak h$ are the Lie algebras of $G$ and
$H$ respectively.

\medskip

Moreover, the map $\ind^G_H$ is uniquely determined by properties (c) and (d).
\end{prop}

\begin{proof}
Naturality of the map $\ind^G_H$ follows from naturality of the action
of the Weyl group described in Section \ref{ss.subgroups}.  
Assertion (b) is a straightforward calculation.  For the
projection formula, if $\alpha \in R(H)$ and $\beta \in G(G,X)$, then
$\alpha \res(\beta)$ is identified with the element $\alpha \beta \in
G(T,X)$.  Since $\beta$ is $W$-invariant, $w(\alpha \beta) = w(\alpha)
\beta$; the projection formula follows.  We now turn to the
reciprocity formula.  Using the projection formula we reduce to the
case where $\beta$ is trivial; then we want to show
$$ (\ind^G_H(\alpha))^G = (\lambda_{-1}({\mathfrak g^*}/{\mathfrak
h^*}) \alpha)^H .$$ 
First assume that $H = T$ so by Serre's theorem we may identify 
$R(G)$ with $R(T)^W$.  
The groups $H$ and $G$ are complexifications of
their corresponding maximal compact Lie subgroups \cite[Theorem 8, p. 244]{OnVi:90},
and
a compact Lie group has the same representation ring as its
complexification \cite[Ch. III, Prop.~8.6]{BrtD:95}.
The Weyl integration formula \cite[Ch. IV, Theorem 1.11]{BrtD:95} for compact Lie groups implies that if
$\nu \in R(G)$, then
\begin{equation} \label{e.recip1}
\nu^G = \frac{1}{|W|}( (\lambda_{-1}({\mathfrak g^*}/{\mathfrak t^*}) \nu)^T.
\end{equation}
If $\nu = \ind_T^G(\alpha)$, then, since $\lambda_{-1}({\mathfrak
g^*}/{\mathfrak t^*})$ is $W$-invariant, the right hand side of
\eqref{e.recip1} becomes
$$
\frac{1}{|W|} \sum_{w \in W} ( w (\lambda_{-1}({\mathfrak g^*}/{\mathfrak t^*}) \alpha))^T.
$$ 
But if $r \in R(T)$ then $r^T = (wr)^T$ for all $w \in W$; the
reciprocity formula follows for $H=T$.  For general $H$, we suppose $T
\subset H \subset G$.  We may assume $\alpha = \ind_T^H(\beta)$.
Then, applying (b) and what we have just proved,
$$ (\ind_H^G(\alpha))^G = (\ind_T^G(\beta))^G =
(\lambda_{-1}({\mathfrak g^*}/{\mathfrak t^*}) \alpha)^T.
$$
On the other hand, applying the projection formula and what we have just proved with $H$ in place of $G$, we obtain
$$ (\lambda_{-1}({\mathfrak g^*}/{\mathfrak h^*}) \alpha)^H =
(\ind_T^H(\lambda_{-1}({\mathfrak g^*}/{\mathfrak h^*}) \beta))^H =
(\lambda_{-1}({\mathfrak h^*}/{\mathfrak t^*})\lambda_{-1}({\mathfrak
g^*}/{\mathfrak h^*}) \beta)^T.
$$ Since $\lambda_{-1}({\mathfrak g^*}/{\mathfrak t^*}) =
\lambda_{-1}({\mathfrak h^*}/{\mathfrak t^*}) \lambda_{-1}({\mathfrak
g^*}/{\mathfrak h^*}) $, the reciprocity formula follows.  Finally,
the projection formula and Proposition \ref{prop.tensor} imply that
the map $\ind_H^G$ is uniquely determined by its value when $X$ is a
point, in which case it is determined by the reciprocity formula.
\end{proof}

\subsection{Induction for Levi factors of parabolic subgroups}
If $H$ is a Levi factor of a parabolic subgroup $P$ of $G$, and $X$
is a $G$-space, then
$$
G(H,X) \simeq G(G,G/H \times X) \simeq G(G,G/P \times X).
$$
Here the first isomorphism was explained in Section \ref{ss.subgroups},
and the second is because the projection $G/H \times X \to G/P \times X$
has fibers isomorphic to affine space.  

We include a proof of the following lemma because of the lack of a reference.

\begin{lemma}
Let $G$ be a connected reductive group.  The commutator subgroup $G'$ of $G$ is simply connected
if and only if $\pi_1(G)$ is torsion-free.  
\end{lemma}

\begin{proof}
The group $G$ is the product
of the semisimple group $G'$ and a torus $S$, with $G' \cap S$ finite
(see \cite[Section 0.9]{Hum:95}
Therefore $G/G'$ is isomorphic
to the torus $S/(G' \cap S)$.  Since $\pi_2(G/G') = 0$, the long exact homotopy sequence of the
fibration $G \to G/G'$ yields a short exact sequence
$0 \to \pi_1(G') \to \pi_1(G) \to \pi_1(G/G') \to 0$.  
Since $G/G'$ is a torus, $\pi_1(G/G')$ is free abelian.
Since $G'$ is semisimple, $\pi_1(G')$ is finite \cite[Chapter VII, Theorem
6.1]{Hel:78}.
The lemma follows.
\end{proof}

\begin{prop} \label{prop.simplelevi}
Let $H$ be a Levi factor of a parabolic subgroup of $G$.  
If $G$ has simply connected commutator subgroup, then so does $H$.
\end{prop}

\begin{proof}
 The quotient $G/H$ has the homotopy type of $G/P$,
so $G/H$ is simply connected and by the Hurewicz theorem,
$\pi_2(G/H)$ is isomorphic to $H_2(G/H) \simeq H_2(G/P)$, which
is torsion-free.  Since $\pi_2(G) = 0$,
the long exact homotopy sequence of the fibration $G \to G/H$
implies that $\pi_1(H)$ is torsion-free.
\end{proof}

In this situation,
the following proposition gives another construction of the twisted
induction map.
\begin{prop} \label{p.altinduction}
Let $G$ and $H$ be as in Proposition \ref{prop.simplelevi}.  Identify $G(H,X)$ with $G(G,G/P \times X)$,
and let $q: G/P \times X \to X$ denote the projection.  If $\alpha \in G(H,X)$,
then
\begin{equation} \label{e.alternateind}
\ind_H^G(\alpha) = q_*(\lambda_{-1}({\mathfrak g^*}/{\mathfrak p^*}) \alpha).
\end{equation}
\end{prop}

\begin{proof}  
The map defined by the right hand side of \eqref{e.alternateind}
satisfies the projection formula.  Hence by Propositions \ref{prop.tensor} and 
\ref{prop.weylinduction} it suffices to show \eqref{e.alternateind} hold for
$\alpha \in R(H)$.  When computing the right hand side of
\eqref{e.alternateind} we view $\lambda_{-1}({\mathfrak
g^*}/{\mathfrak p^*}) \alpha$ as inducing a virtual vector bundle on
$G/P$, and $q_*$ as taking the $G$-equivariant Euler characteristic of
that bundle.  Here, if $V$ is a representation of $H$, then we may
view it as a $P$-module by making the unipotent radical act trivially.
The induced vector bundle on $G/P$ is $G \times^P V$.  The vector
bundle on $G/P$ induced by ${\mathfrak g^*}/{\mathfrak p^*}$ is the
cotangent bundle $T^* = T^*_{G/P}$, so what we need to compute is
$q_*(\lambda_{-1}(T^*) \alpha$).  The latter computation may be done
using the same technique as in the proof Proposition 3.10 of
\cite{EdGr:05}, as we now explain.

The restriction
map $G_0(G,G/P) \to G_0(T, G/P)$ is injective,
so we may restrict to $T$-equivariant $K$-theory, and make the calculation in
$G_0(T,G/P)$.   In particular if $a \in T$ is any element
such that $(G/P)^a = (G/P)^T$, then the localization theorem in $T$-equivariant
$K$-theory \cite[Theorem 3.3]{EdGr:05} implies that
\begin{equation}  \label{e.alternateind2}
\lambda_{-1}(T^*) \alpha = \sum_{x \in (G/P)^T} i_{x*} \frac{i_x^*(\lambda_{-1}(T^*) \alpha)}{\lambda_{-1}(T_x^*)} = \sum_{x \in (G/P)^T} i_{x*} i_x^*{\alpha}
\end{equation}
in $G_0(T,G/P)_{{\mathfrak m}_a}$,
where $i_x$ is the inclusion of $x$ in $G/P$.
Since $G/P$ has a cellular decomposition, 
\cite[Lemma 5.5.1]{ChGi:97} implies that $G_0(T,G/P)$ is a free $R(T)$-module
of rank equal to the number of $T$-fixed points. Since $R(T)$ is an integral
domain, it follows that the equality
\begin{equation}\label{e.referee}
\lambda_{-1}(T^*) \alpha = \sum_{x \in (G/P)^T} i_{x*} i_x^*{\alpha}
\end{equation}
also holds in $G_0(T,G/P)$.
The $T$-fixed points in $G/P$ are the points $wP$, for $w \in W/W_1$.
If $x$ corresponds to the coset $wP$ then, 
identifying $G(T, x)$ with $R(T)$, the map $q_* \circ  i_{x}$ is the identity,
and for $\alpha \in R(H) = R(T)^{W_1}$, we have $i_{x}^*(\alpha) = w \alpha$.
Applying $q_*$ to both sides of \eqref{e.referee} yields
$$
q_*(\lambda_{-1}(T^*) \alpha) = \sum_{w \in W/W_1} w \alpha = \ind_H^G{\alpha},
$$
as desired.
\end{proof}

\section{Induction, restriction and completion}
Let $G$ be a connected reductive group with simply connected commutator subgroup
and let $\Psi \subset G$ be a semisimple 
conjugacy class. Let $h$ be an element of $\Psi$, and let $Z = {\mathcal Z}_G(h)$ be the
centralizer of $h$.  Assume that $Z$ is also connected and
reductive and has simply connected commutator subgroup.
\footnote{If $G$ has simply connected commutator subgroup then $Z$ is automatically connected (c.f.
\cite[Theorem~2.11]{Hum:95}).} In this
section we prove that, after appropriate completions, the
induction functor $\ind^G_Z$ is in fact the inverse of the restriction
functor.

\subsection{Some facts about completions and invariants}
In this section we prove some basic commutative algebra results about invariants
and completions which will be necessary for proving properties of the 
induction map in equivariant  $K$-theory.

\subsubsection{Completions and finite morphisms}
Let $B \to A$ be a finite morphism of Noetherian rings. Let ${\mathfrak
m}$ be a maximal ideal of $B$, let ${\mathfrak a}= \mathfrak{m}A$,
and let ${\mathfrak m}_1, \ldots
\mathfrak{m}_r$ be the maximal ideals of $A$ lying over ${\mathfrak
m}$. Let $M$ be an $A$-module and let $\widehat{M}$ be the ${\mathfrak
a}$-adic completion of $M$.
Likewise, let $\widehat{M}_i$ denote the ${\mathfrak m}_i$-adic completion
of $M$. Since ${\mathfrak m}M \subset {\mathfrak m_i}M$ there is an
induced map of completions $v_i \colon \widehat{M} \to \widehat{M}_{i}$.  The
following proposition will be needed. It is a generalization of
\cite[Chapter III, Section 2, no. 13, Proposition 18]{Bourbaki}.
\begin{lemma} \label{lem.completions}
The map$$ (v_i) :  \widehat{M} \to \prod_{i = 1}^r \widehat{M}_{i}$$
is an isomorphism of $\widehat{A}$-modules.
\end{lemma}

\begin{proof}
Let $N$ denote the product $\prod_{i =1}^r M$ with the topology given by
the product of the ${\mathfrak m}_i$-adic topologies, $i =1 ,\ldots r$.
The completion of $N$ in this topology is the product
$\prod_{i =1}^r \widehat{M}_i$.

We claim that the
${\mathfrak a}$-adic topology on $M$ is the same as the topology
induced by the inclusion by the diagonal map
$M \stackrel{\Delta} \to N$. To see this argue
as follows.  Let $N_k = \prod_{i =1}^r {\mathfrak m_i}^k M \subset N$.
The filtration $N \supset N_1 \supset \ldots \supset N_k \supset
\ldots$ is a fundamental system of neighborhoods of $0$ so the
completion of $N$ can be viewed as the completion with respect to this
filtration.  Now $\Delta^{-1}(N_k) = \cap {\mathfrak m_i}^k M =
{\mathfrak m_1}^k \ldots {\mathfrak m_r}^k M$ so the induced topology 
on $M$ is the ${\mathfrak m}_1 \ldots {\mathfrak m}_r$ -adic
topology. But $\sqrt{{\mathfrak a}} = {\mathfrak m}_1 \ldots {\mathfrak m}_r$ so
this coincides with the ${\mathfrak a}$-adic topology as $A$ is
Noetherian.

If $n= (m_1, \ldots , m_r) \in N$ then by the Chinese remainder theorem
for modules
\cite[Chapter II, Section 1, no. 2, Proposition 6]{Bourbaki}
there exists $m \in M$ such that $m \equiv m_i \mod {\mathfrak m_i}^nM$
so $n \in \Delta(M) + N_n$ for any $n$. Hence $\Delta(M)$ is dense in $N$.

Let $N/\Delta(M)$ have the topology induced from the topology on
$N$. Since the ${\mathfrak a}$-adic topology on $M$ is induced
from the topology on $N$, the sequence of completions
$$0 \to \widehat{M} \stackrel{\widehat{\Delta}} \to \widehat{N} \to
\widehat{N/\Delta(M)} \to 0$$ is exact \cite[Chapter III, Section 12,
no. 12 Lemma 2]{Bourbaki}. But $\Delta(M)$ is dense, so
$\widehat{N/\Delta(M)} = 0$. Identifying the completion $\widehat{N}$ with
$\prod_{i = 1}^r \widehat{M_i}$ and $\widehat{\Delta}$ with $\prod_{i =
1}^r{v_i}$ yields the proposition.
\end{proof}

\subsubsection{Completions and invariants}
Let $W$ be a finite group acting on a ring $A$.  
Assume that $|W|$ is invertible in $A$. Let
$B = A^W$  denote the subring of invariants.  Let $I \subset B$
be an ideal, and $\widehat{B}$ the $I$-adic completion of $B$.

Let $M$ be an $A$-module with a $W$-action compatible with
the module structure, and let $M^W$ be the $B$-submodule
of  $W$-invariants.
Since $I \subset B$ is $W$-invariant, there is an action
of $W$ on the $I$-adic completion
$\widehat{M}$ of $M$.
\begin{lemma} \label{lem.technicalcompletions}
Let $\widehat{M^W}$ be the $I$-adic completion of $M^W$.
There is a natural isomorphism of $\widehat{B}$-modules
$\widehat{M^W} \to (\widehat{M})^W$.
\end{lemma}
\begin{proof}
First observe that since $W$ acts trivially on $\widehat{B}$,
$(\widehat{M})^W$ is in fact a $\widehat{B}$-module.  The hypothesis that $1/|W| \in A$ allows us to define a projection
$M \to M^W$ by the formula $m \mapsto {1\over |W|}\sum_{w \in W} wm$.
Thus,
for any integer $k >0$, the exact sequence of $A$-modules
$$0 \to I^k M \to M \to M/I^kM \to 0$$
induces an exact sequence of $B$-modules
$$0 \to (I^kM)^W \to M^W \to (M/I^kM)^W \to 0.$$
Moreover, $(I^kM)^W = I^kM^W$, so the exact
sequence induces isomorphisms
\begin{equation} \label{e.technical}
 M^W/I^kM^W \to (M/I^kM)^W
 \end{equation}
 for all $k$.
 By definition,
 $$
 \liminv  M^W/I^kM^W = \widehat{M^W}.
 $$
Also, an element $(m_k)\in \widehat{M}$ is $W$-invariant 
if and only if each $m_k$ is $W$-invariant, i.e.,  
$$(\widehat{M})^W = \liminv (M/I^kM)^W
\subset \liminv M/I^kM = \widehat{M}.$$
Thus the desired isomorphism $\widehat{M^W} \to (\widehat{M})^W$
follows from \eqref{e.technical} by taking inverse limits.
\end{proof}

\subsubsection{Invariants and direct products} \label{sec.invdirect}
Let $A$ be a ring, let $M_1  \ldots, M_r$ be $A$-modules,
and let $M = \prod_{i =1}^r M_i$.
Suppose that a finite group $W$ acts compatibly on $A$ and $M$.  Assume
that for all $w \in W$, $w M_i = M_j$ for some $j$, and conversely
that for any $i$ and $j$ there exists $w \in W$ such that $w M_i =
M_j$.  Let $W_i = \{w \in W \ | \ wM_i = M_i\}$. The transitivity of the
action on the components of $M$ implies that the $W_i$ are conjugate
subgroups of $W$.  Let $\pi_i \colon M \to M_i$ denote the projection.
If $m \in M^W$ then the $i$-th component of $m$ must be in
$M_i^{W_i}$, so $\pi_i$ restricts to a map $\res_i: M^W \to
M_i^{W_i}$.

Let $w_1, \ldots w_r$ be a set of representatives for the
cosets of $W/W_i$ such that $w_j M_i = M_j$.   We define
a map $\ind_i:M_i^{W_i} \to M^W$ by  $\ind_i(m_i) \mapsto (w_1m_i, \ldots
w_rm_i)$.
\begin{lemma} \label{lem.invdecomp}
With the above assumptions and notation,
the maps $\res_i$ and $\ind_i$ are inverse isomorphisms of
the $A^W$-modules $M^W$ and $M_i^{W_i}$.
\end{lemma}
\begin{proof}
Let $m \in M^W$. Since $W$ acts transitively on the components
of $M$, $\res_i(m) = 0$ if and only if $m = 0$, so $\res_i$ is injective.
Since $\res_i \circ \ind_i$ is the identity, $\res_i$ and $\ind_i$ are inverse
isomorphisms.
\end{proof}
\begin{remark}
Our hypothesis on the $W$ action on $M$ implies 
that $M = \Ind^W_{W_i}(M_i)$ for any $i$, where $\Ind$ is the standard
induction functor from $W_i$-modules to $W$-modules.  Viewed this way,  Lemma \ref{lem.invdecomp}
is simply the statement that $(\Ind^W_{W_i}(M_i))^{W} = M_i^{W_i}$.
\end{remark}
\subsection{Induction and restriction} \label{sec.indres}
Let $G$ be a complex algebraic group and $\Psi \subset G$ a semisimple
conjugacy class. Fix $h \in \Psi$ and set $Z = {\mathcal Z}_G(h)$. Let
$\m_\Psi \subset R(G)$ be the maximal ideal of virtual representations
whose character vanishes on $\Psi$ and $\m_h \subset R(Z)$ the maximal
ideal of virtual representations whose character vanishes at the
central element $h$.  By \cite[Propositions 2.3, 2.6]{EdGr:05} the
restriction map $R(G) \to R(Z)$ is finite and $\m_h$ lies over
$\m_\Psi$.  If $M$ is an $R(G)$-module, let $\widehat{M}$ denote the
$\m_\Psi$-adic completion of $M$.  If $M$ is also an $R(Z)$-module,
such that the $R(G)$-module structure is obtained from the
$R(Z)$-module structure by the map $R(G) \to R(Z)$, then Lemma
\ref{lem.completions} implies that the $\m_h$-adic completion
$\widehat{M}_h$ of $M$ is canonically identified as a summand in the
$\m_\Psi$-adic completion $\widehat{M}$. The natural map $\widehat{M}
\to \widehat{M}_h$ corresponds to projection onto this summand.

If $X$ is a $G$-space then, in the context of equivariant $K$-theory, 
we have maps
$$\res_{h} \colon \widehat{G(G,X)} \to
\widehat{G(Z,X)} \to \widehat{G(Z,X)}_{h}
$$
and
$$\res_h \colon \widehat{K_0(G,X)} \to \widehat{K_0(Z,X)} \to
\widehat{K_0(Z,X)}_{h}$$ corresponding to the composition of the
restriction maps $\widehat{G(G,X)} \to \widehat{G(Z,X)}$
(resp. $\widehat{K_0(G,X)} \to \widehat{K_0(Z,X)}$) with projection on
the summand $\widehat{G(Z,X)}_h$ (resp. $\widehat{K_0(Z,X)}_h$).  If
in addition we assume that both $G$ and $Z$ are connected and
reductive and have simply connected commutator subgroups, then the
induction map $\ind: G(Z,X) \to G(G,X)$ induces a map on completions
at the ideal $\m_{\Psi}$.  Composing with the inclusion of
$\widehat{G(Z,X)}_h$ as a summand in $\widehat{G(Z,X)}$ gives a map $$
\ind_{h}: \widehat{G(Z,X)}_{h} \to \widehat{G(G,X)}.
$$

\begin{prop} \label{prop.indres}
If $G$ and $Z$ are connected, reductive and have simply connected commutator
subgroups, then
the maps $\res_h$ and $\ind_h$ are inverse isomorphisms of
$\widehat{R(G)}$-modules.
\end{prop}

Before we prove Proposition \ref{prop.indres} we need a lemma about
Weyl groups of centralizers of semisimple elements.
\begin{lemma}\label{lem.weyl}
Let $G$ be a connected reductive group and let $h$ be a semisimple
element of $G$ such that the centralizer $Z = {\mathcal Z}_G(h)$ is
connected.  Let $\Psi$ denote the conjugacy class of $h$ in $G$.  Let
$T$ be a maximal torus of $Z$ (hence also of $G$).  Let $W =
W(G,T)$and $W_1 = W(Z,T)$.  Then $h \in T$, $\Psi \cap T = W \cdot h$, and the
stabilizer in $W$ of $h$ is $W_1$.
\end{lemma}

\begin{proof}[Proof of Lemma \ref{lem.weyl}]
First, $h$ is
contained in some maximal torus of $Z$.  The conjugacy
of maximal tori in $Z$ implies that the the $Z$-conjugacy class
of $h$ meets $T$, but this class consists of $\{ h \}$.  Hence
$h \in T$.  Suppose that $h' = ghg^{-1} \in T$.  Since $g^{-1}Tg = T_1$ and $T$
are two maximal tori containing $h$, they are both contained in $Z$.
Therefore there is an element $z \in Z$ such that $zTz^{-1} = T_1$.
Then $h' = gzh(gz)^{-1}$.  The element $gz$ of $N(T)$ represents an
element of $W = N(T)/T$, so $h' \in W \cdot h$.  Hence $\Psi \cap T
\subseteq W \cdot h$; the reverse inclusion holds as well since
elements of $W$ are represented by elements of $N(T) \subset G$.  This
proves the first statement.  For the second, suppose $w \in W$ fixes
$h$.  Then $w$ is represented by an element of $Z$, so $w$ lies in the
Weyl group of $Z$, namely $W_1$.
\end{proof}

\begin{proof}[Proof of Proposition \ref{prop.indres}]
Choose a maximal torus $T$ of $Z$ (and $G$). Let $W = W(G,T)$ and
$W_1 = W(Z,T)$.  By Lemma \ref{lem.weyl},
the Weyl group $W$ acts transitively on the set
$\Psi \cap T$ and the stabilizer of $h$ is $W_1$.  
Thus, $\Psi \cap T = \{h_1, h_2, \ldots h_r\}$ where $r =
{|W| \over{ |W_1|}}$ (we take $h=h_1$).   
Set $M = G(T,X)$ and let $\widehat{M}_i$ be the completion of $M$ at
the maximal ideal
${\mathfrak m}_{h_i}$ of $R(T)$.  By Lemma \ref{lem.completions},
$\widehat{M} \simeq \oplus_{i=1}^r \widehat{M}_i$.
By Lemma \ref{lem.technicalcompletions},
$\widehat{G(G,X)} = (\widehat{M})^W$ and $\widehat{G(Z,X)}_h =
(\widehat{M}_1)^{W_1}$.  The map $\ind_h$ corresponds to
the induction map $(\widehat{M_1})^{W_1} \to (\widehat{M})^W$ defined in Section
\ref{sec.invdirect} above, and the map $\res_h$ corresponds to
the projection $(\widehat{M})^W \to (\widehat{M}_1)^{W_1}$.  Hence Lemma
\ref{lem.invdecomp} implies that $\ind_h$ and $\res_h$ are inverse isomorphisms.
\end{proof}

\begin{example} \label{ex.needcompletions}
Proposition \ref{prop.indres} illustrates the need to work with completions
rather than localizations, because the restriction map of localizations
$G(G,X)_{\m_\Psi} \to G(Z,X)_{\m_h}$ will not generally be an isomorphism.
For example, if $G = \GL_n$ and $\Psi$ is a regular conjugacy class (i.e.
elements of $\Psi$ have distinct eigenvalues), then $Z = T$ is a maximal
torus.  The restriction map $R(\GL_n) \to R(T)$ has degree $n!$, so the
map of local rings $R(G)_{\m_\Psi} \to R(T)_{\m_h}$
cannot be an isomorphism. However, this map of local rings is \'etale
so the corresponding map of completions is indeed an isomorphism.
\end{example}

\section{The non-abelian completion theorem}
In this section we state and prove Theorem \ref{thm.nonabelian}, our
general non-abelian completion theorem for equivariant
$K$-theory. This theorem extends the non-abelian localization theorem
of \cite[Theorem 5.1 and Corollary 5.2]{EdGr:05} which was proved for
actions with finite stabilizer.

\subsection{The completion theorem for a central element}
Let $Z$ be an algebraic group and let $h$ be an element
in the center of $Z$. If $M$ is an $R(Z)$-module, we denote by
$\widehat{M}_h$ the ${\mathfrak m}_h$-adic completion of
$M$, where, as usual, ${\mathfrak m}_h$ is the ideal in $R(H)$
corresponding to virtual representations whose characters vanish at
$h$.

If $X$ is a $Z$-space, then
since $h$ is central, the fixed locus $X^h$ is $Z$ invariant.
When $X$ is smooth, then so is $X^h$, and thus the inclusion
$i_* \colon X^h \to X$ is a regular embedding.
As a tool in proving our general completion theorem, we need the following
result.

\begin{thm} \label{thm.centralcompletion}
Let $Z$ be an algebraic group and $h$ a central element of $Z$.
Let $X$ be a $Z$-space and let $i \colon X^h \to X$ denote
the inclusion of the fixed point locus of $h$.

(a) The proper pushforward $i_* \colon G(Z,X^h) \to G(Z,X)$
is an isomorphism after completing at ${\mathfrak m}_h$.

(b) If $X$ is smooth, then the map of $R(Z)$-modules
$$\cap \; \lambda_{-1}(N_i^*) \colon G(Z, X^h) \to G(Z,X^h)$$
is invertible after completing at ${\mathfrak m}_h$, and
if $\alpha \in \widehat{G(Z,X)}_h$ then
\begin{equation} \label{eq.abcompletion}
\alpha = i_*(\lambda_{-1}(N_i^*)^{-1} \cap i^*\alpha),
\end{equation}
where the notation $\lambda_{-1}(N_i^*)^{-1} \cap i^*\alpha$
means inverse image of $i^*\alpha$ under the isomorphism
$\cap \; \lambda_{-1}(N_i^*)$.
\end{thm}
\begin{proof}
The proof of Theorem \ref{thm.centralcompletion}
is essentially the same as the proof of the central localization theorem
(\cite[Theorem 3.3]{EdGr:05}). Details may be found in that paper; here is
a brief sketch.
By a change of groups argument the result for a general group
reduces to the case of a product of general linear groups, and
this in turn reduces to the case of a maximal torus.
Therefore we may assume that $G=T$ is a torus.
Using
the localization long exact sequence for an open set,
to prove part (a) it suffices to show that if $X^h$ is empty
then $\widehat{G(T,X)}_h = 0$. By Thomason's generic slice theorem
we may reduce to the case that $X= T/T' \times Y$ where $T' \subset
T$ is a closed subgroup not containing $h$, and the action of $T$ on $Y$ is trivial.
In this case $G(T,X) = R(T') \otimes G(Y)$ and $\widehat{G(T,X)}_h
= \widehat{R(T')}_h \otimes G(Y)$. The crucial point is that
if $h \notin T'$ then $R(T')_{{\mathfrak m}_h} = 0$ (which gives the
localization result)
and $\widehat{R(T')}_h = 0$ (which gives the completion result).
This proves (a).
To prove (b), we again use the generic
slice theorem to reduce to the case where $X^h = T/T' \times X^h/T$, but now $T'
\subset T$ is a closed subgroup containing $h$. Then $G(T,X^h) = R(T')
\otimes G(X^h/T)$ so $\widehat{G(T,X^h)}_h = \widehat{R(T')}_h \otimes
G(X^h/T)$. Let ${\mathcal N}_x$ denote the fiber of $N_i$ at a point
$x \in X^h$.
The action of $\lambda_{-1}(N_i^*)$ on $\widehat{G(T,X^h)}_h$
is invertible if and only if $\lambda_{-1}({\mathcal N}_x)$ is invertible in
$\widehat{R(T)}_h$ for some closed point $x$ in each connected component
of $X^h$. This follows from the fact (proved in \cite{EdGr:05}) that
$\lambda_{-1}({\mathcal N}_i) \notin {\mathfrak m}_h$.
Finally, to deduce \eqref{eq.abcompletion}, if $\alpha \in \widehat{G(Z,X)}_h$,
then by (a) we can write $\alpha = i_* \beta$ for 
$\beta \in \widehat{G(Z,X^h)}_h$.  By the equivariant self-intersection
formula, $\lambda_{-1}(N_i^*) \beta = i^* \alpha$, so
$\beta = \lambda_{-1}(N_i^*)^{-1} \cap i^* \alpha$.
\end{proof}

\begin{remark}
In the proof of the central localization theorem of \cite{EdGr:05},
the change of groups step proceeded by embedding
$G$ into a product of general linear groups $Q$ such that
the conjugacy class of $h$ in $Q$ intersects $G$ in the single
point $h$.  We asserted that $h$ central in $G$ implies $h$ central in
$Q$.  Unfortunately this assertion can fail for the embeddings
we constructed in that paper.  
However, we can arrange it to be true if $G$ is
reductive, since then we can take $Q$ to be $\prod GL(V_i)$, where 
$V_i$ is an irreducible $G$-module; in this case the assertion
follows from Schur's lemma, and the proof of the central
localization theorem goes through.  To deal with nonreductive $G$
we replace $G$ by a Levi factor $L$ and use the identification of
$G(G,X)$ with $G(L,X)$; then the argument works because
$L$ is reductive.
\end{remark}

\subsection{Statement of the non-abelian completion theorem}
Let $X$ be a $G$-space
and $\Psi$ a semisimple conjugacy class in $G$.  Let $h$ be an
element of $\Psi$ and $Z = {\mathcal Z}_{G}(h)$.
As usual, let $i \colon X^h \to X$ be the inclusion of the fixed point locus of $h$;
as noted above,
if $X$ is smooth then so is
$X^h$, and the inclusion is a regular embedding.
Since $Z$ acts on $X^h$, if $X$ is smooth we may define
a map $i^! \colon \widehat{G(G,X)} \to \widehat{G(Z,X^h)}_{h}$
as the composition
$$
\widehat{G(G,X)} \stackrel{\res_h}{\rightarrow} \widehat{G(Z,X)}_h
\stackrel{i^*}{\rightarrow} \widehat{G(Z,X^h)}_{h},
$$
where $\res_h$ is as in Section \ref{sec.indres}.

For arbitrary $X$ we may also define $i^! \colon \widehat{K_0(G,X)}
\to \widehat{K_0(Z,X^h)}_{h}$ in the same way.  This pullback will be
used in the statement of Theorem \ref{thm.rrnextgeneration}.

\begin{thm} \label{thm.nonabelian}
Let $X$ be a $G$-space, and use the notation above.
There is an
isomorphism $i_! \colon \widehat{G(Z,X^h)}_{h}
\to \widehat{G(G,X)}$ with the following
properties.

(a) (covariance) Let $Y \stackrel{p} \to X$ be a proper $G$-morphism
of algebraic spaces, and let $Y^h \stackrel{q} \to X^h$ be the
corresponding proper map of $h$-fixed loci. Then
$ i_! q_* = p_* i_!$ as maps $\widehat{G(Z,Y^h)}_{h} \to \widehat{G(G,X)}$.

(b) If $X$ is smooth, then $i^!$ is also an isomorphism
and we have the following analogue of the explicit localization formula:
\begin{equation} \label{eqn.nonabelian}
\alpha = i_!
\left( \lambda_{-1}(N_i^*)^{-1} \cap i^!\alpha\right)
\end{equation}

(c) The map $i_! \colon \widehat{G_0(Z,X)}_{h} \to \widehat{G_0(G,X)}$
is uniquely determined by by properties (a) and (b).
\end{thm}

\begin{remark}
We conjecture that the uniqueness statement of part (c) holds for
higher $K$-theory as well; this would follow if we had a Riemann-Roch theorem 
relating higher equivariant $G$-theory
to higher equivariant Chow groups. 
\end{remark} 
When $G$ and $Z$ are connected and reductive with simply connected commutator
subgroup, then Proposition 
\ref{prop.indres} states that for any $G$-space $X$, the map
$\res_h \colon \widehat{G(G,X)} \to \widehat{G(Z,X)}_h$ is an isomorphism. Using Theorem \ref{thm.nonabelian},
we can partially generalize this to arbitrary $G$ and $Z$:

\begin{cor} \label{cor.nonabelian}
(a) If $X$ is a smooth $G$-space, then 
$\res_h \colon \widehat{G(G,X)} \to \widehat{G(Z,X)}_h$
is an isomorphism. In particular, the restriction map
$R(G) \to R(Z)$ induces an isomorphism between the 
${\mathfrak m}_\Psi$-adic completion of $R(G)$ and the
${\mathfrak m}_h$-adic completion of $R(Z)$.

(b) For any $G$-space $X$, the restriction map of completed
equivariant Grothendieck groups 
$\res_h \colon \widehat{G_0(G,X)} \to \widehat{G_0(Z,X)}_h$ is an isomorphism.
\end{cor}
\begin{proof}
If $X$ is smooth, then the central completion theorem implies that
$i^* \colon \widehat{G(Z,X)}_h \to \widehat{G(Z,X^h)}_h$
is an isomorphism. Since $i_!$ is an isomorphism, 
formula \eqref{eqn.nonabelian} of Theorem \ref{thm.nonabelian} implies that
$i^!$ is an isomorphism. Hence $\res_h$ must also be an isomorphism,
proving (a).  The proof of (b)
uses envelopes and Riemann-Roch,
and will be deferred to Section \ref{sec.envelopeproof}.
\end{proof}

We conjecture that the conclusion of part (a) holds for arbitrary $X$.

\subsection{Proof of Theorem \ref{thm.nonabelian}(a),(b)}

In this section we prove parts (a) and (b) of Theorem \ref{thm.nonabelian}.
Part (c) uses a Chow envelopes argument and
our generalized equivariant Riemann-Roch theorem;
it will be proved in Section \ref{sec.envelopeproof}.

If $G$ and $Z$ are connected and reductive with simply connected commutator 
subgroups, we define
$i_!$ as $\ind_h \circ i_*$. By the central completion theorem 
(Theorem \ref{thm.centralcompletion})
$i_*$ is an isomorphism, and by Proposition \ref{prop.indres},
$\ind_h$ is an isomorphism. If $Y \to X$ is proper, then the diagram
$$\begin{array}{ccc}
Y^h & \to & Y\\\downarrow & & \downarrow\\
X^h & \rightarrow & X
\end{array}
$$
is a commutative diagram of proper morphisms. Covariance of $i_!$
then follows from covariance of both pushforwards and $\ind_h$
for proper $G$-morphisms.  Part (b) follows from Theorem \ref{thm.centralcompletion}(b)
and the fact
that $\ind_h \circ \res_h = id$.

To prove the general case we use a change of groups argument similar
to that used in \cite{EdGr:05}. By \cite[Proposition 2.8]{EdGr:05}, we
may embed $G$ into a product of general linear groups $Q$ such that,
writing $\Psi = C_G(h)$ and $\Psi_Q = C_Q(h)$, we have $\Psi_Q \cap G
= \Psi$.   

Write $Z= {\mathcal Z}_G(h)$ and $Z_Q = {\mathcal
Z}_Q(h)$. 
The groups $Q$ and $Z_Q$ are both connected and reductive and have 
simply connected commutator
subgroups because they are both products of 
general linear groups.  Let $Y = Q \times_G X$, and let
$j:Y^h \to Y$ be the inclusion of the fixed point locus.
By what we have proved,
there is an isomorphism of completions $j_! \colon \widehat{G(Z_Q,
Y^h)}_h \to \widehat{G(Q,Y)}$.  
Here $\widehat{G(Z_Q,
Y^h)}_h$ denotes the completion of $G(Z_Q,
Y^h)$ at the maximal ideal $\m_h'$ of $R(Z_Q)$
corresponding to the central element $h$, and
$\widehat{G(Q,Y)}$ is the $\m_{\Psi_Q}$-adic completion of
$G(Q,Y)$.

We want to use $j_!$ to define $i_!$.  First observe that the only
maximal ideal of $R(G)$ containing $\m_{\Psi_Q}R(G)$ is $\m_{\Psi}$
(\cite[Proposition 2.6]{EdGr:05}), so $\m_{\Psi_Q}R(G)$ is an
$\m_{\Psi}$-primary ideal.  Thus, under the Morita equivalence which
identifies $G(Q,Y)$ with $G(G,X)$, the ${\mathfrak m}_{\Psi_Q}$-adic
and ${\mathfrak m}_{\Psi}$-adic topologies coincide. Hence we may
identify the ${\mathfrak m}_{\Psi_Q}$-adic completion of $G(Q,Y)$ with
the ${\mathfrak m}_{\Psi}$-adic completion of $G(G,X)$.  Similarly,
under the Morita equivalence identification of $G(Z,X^h)$ with $G(Z_Q,
Z_Q \times_Z X^h)$, the $\m_h$ and $\m_h'$-adic topologies coincide.

The following lemma is a consequence of
an extension of \cite[Lemma 5.5]{EdGr:05} to
non-smooth spaces.

\begin{lemma} \label{lem.changeofgroups}
The natural map $Z_Q \times_Z X^h \to
Y^h$ is an isomorphism of $Z_Q$-spaces.
\end{lemma}
\begin{proof}[Proof of Lemma \ref{lem.changeofgroups}]
Let $S =\{(g,x) \ | \ g\in \Psi , gx =x\}
\subset G \times X$ and
$S_Q = \{(q,y) \ | \ q \in \Psi_Q, qy = y\} \subset Q \times Y$.
As in the proof of \cite[Lemma 5.5]{EdGr:05}, consider the map
$$T \colon Q \times S \to \Psi_Q \times Q \times X,\;\;\; (q,g,x) \mapsto
(qgq^{-1},q,x).$$ This map induces a map of
quotients $\tilde{\Phi} \colon Q \times_G S \to \Psi_Q \times Y$ which factors
through a map $\Phi \colon Q \times_G S \to S_Q \subset \Psi_Q \times Y$.
In view of the identifications $S = G \times_Z X^h$ and
$S_Q = Q \times_{Z_Q} Y^h$ (\cite[Lemma 4.3]{EdGr:05}), to prove the
lemma, it suffices to show that $\Phi:Q \times_G S \to S_Q$ is an isomorphism of $Q$-spaces.
This was proved for smooth $X$ in \cite[Lemma 5.5]{EdGr:05}; to
extend to arbitrary $X$ we argue as follows.  
To check that $\Phi$ is an isomorphism, we may work locally in the
smooth topology on $X$ and assume that $X$ is affine.
\footnote{Every $G$-space
has a smooth $G$-cover by an affine $G$-scheme. To see this, let
$U \stackrel{\pi} \to X$ be any \'etale cover of $X$ by a scheme. Replacing
$U$ by an affine Zariski open cover, we may assume
that $U$ is also affine. Then the map
$G \times U \to X$, $(g,u) \mapsto g\pi(u)$, is a smooth, affine  $G$-cover of
$X$.}
By local linearizability of group actions \cite{Bor:91},
there is a $G$-equivariant embedding $X \subset V$ for some finite
dimensional representation $V$ of $G$. Let $W = Q \times_G V$. Define
$S' = \{(g,v) \ | \ g\in \Psi, gv = v\} \subset \Psi \times V$ 
and $S'_Q = \{(q,w) \ | \ q \in \Psi_Q, qw = w\} \subset Q \times W$. 
Since $V$ is smooth, the map 
$\Phi' \colon Q \times_G S' \to S'_Q$ is an isomorphism.  
Now,
$X^h = X \cap V^h$, so $S = S' \cap (\Psi \times X)$. Hence
the isomorphism $\Phi'$ restricts to the map
$\Phi \colon Q \times_G S \to S_Q$. Therefore, $\Phi$ is an isomorphism.
\end{proof}

Hence we may identify the
${\mathfrak m}_h$-adic completion of $G(Z,X^h)$
with the ${\mathfrak m}_h'$-adic completion of $G(Z_Q, Y^h)$.

Applying Lemma \ref{lem.changeofgroups} to the case $X = \Spec \C$,
we see that $Z_Q/Z$ can be identified with $(Q/G)^h$. Let $e\colon
Z_Q/Z \to Q/G$ be the corresponding $Z_Q$-equivariant regular embedding.
Under the identification of $G(Z_Q,Z_Q/Z)$ with the representation
ring $R(Z)$, the class of
the conormal
bundle $[N^*_e]$ corresponds to the virtual
$Z$-module $\mathfrak{q^*/g^*} - {\mathfrak z_Q}^*/{\mathfrak z}^* = 
{\mathfrak q}^*/{\mathfrak z}_Q^*
- \mathfrak{g^*/z^*}$, where ${\mathfrak q}, {\mathfrak z_Q},{\mathfrak
g}$, and ${\mathfrak z}$ denote the Lie algebras of the groups
$Q,Z_Q,G$, and $Z$, respectively.
The central completion theorem implies that
$\lambda_{-1}(N_e^*)$ is invertible
in $\widehat{G(Z_Q,Z_Q/Z)}_h = \widehat{R(Z)}_h$.

To complete the proof of the theorem for general $G$, we use
the
identifications of $\widehat{G(Z,X^h)}_h$ with $\widehat{G(Z_Q,Y^h)}_h$ and 
$\widehat{G(G,X)}$ with $\widehat{G(Q,Y)}$ to
 define
$$i_! \colon \widehat{G(Z,X^h)}_h \to  \widehat{G(G,X)}$$
by the formula
\begin{equation} \label{e.changeofgroupsstep} 
 i_! \beta = j_!(
\lambda_{-1}(N_{e}^*)^{-1} \cap
\beta  )
\end{equation}
Because $j_!$ is an isomorphism, $i_!$ is as well, and the
covariance of $j_!$ implies covariance if $i_!$.  

Now suppose that $X$ is smooth.
To verify formula \eqref{eqn.nonabelian}
of Theorem \ref{thm.nonabelian},
observe that $N_{j}^* = N^*_i + N_e^*$ in
$G(Z,X^h) = G(Z_Q, Y^h)$, because
$j$ factors as the composition of embeddings
$$
Y^h = (Q \times_G X)^h = (Z_Q \times_Z X^h) \to Z_Q \times_Z X \to
 Y = Q \times_G X.
 $$ 
Moreover, under our identifications, the maps $j^!$ and $i^!$ coincide.  Thus,
formula \eqref{eqn.nonabelian} for $i_!$ follows from the corresponding
formula for $j_!$.

\begin{remark} A priori, our definition of $i_!$ depends
on the choice of the embedding of $G$ into $Q$.
However, if $X$ is smooth, then
 $i_!$ is independent of the embedding, because
equation \eqref{eqn.nonabelian} implies that composition of 
$i_!$ with the map $\cap \;\lambda_{-1}(N_i^*)^{-1}$ is the inverse of $i^!$.
Since $i^!$ is defined intrinsically, so is $i_!$ in this case. 
For arbitrary $X$, the uniqueness
statement Theorem \ref{thm.nonabelian}(c) implies
that $i_!$ is independent of $Q$ as a map of completed
Grothendieck groups.   As noted above, we conjecture that
Theorem \ref{thm.nonabelian}(c) extends to higher $K$-theory;
this would imply that $i_!$ is independent of $Q$ in the higher
$K$-theory case as well.
\end{remark}

\subsection{The case of finite stabilizer}
In this section we show that if $G$ acts with finite stabilizer
on a smooth space $X$ then the nonabelian completion theorem is equivalent to the
nonabelian localization theorem of \cite{EdGr:05}.
In the finite stabilizer case, the equivariant $K$-theory is
supported at a finite number of maximal ideals, so
localization and completion coincide.  Therefore the issue is to compare
the map used in this paper with the
one of \cite{EdGr:05}.

Let $G$ be an algebraic group acting with finite stabilizer
on a smooth algebraic space $X$.
Following the notation of \cite{EdGr:05}, let $\Psi$ be a semisimple conjugacy class in $G$ and 
let $S_\Psi = \{(g,x)|g \in \Psi, gx = x\}$.  
Choose 
$h \in \Psi$, and let $Z =
{\mathcal Z}_G(h)$. 
Then $S_\Psi = G\times_Z X^h$, so by Morita equivalence we may identify $G(G,S_\Psi)$ with 
$G(Z,X^h)$, 
which is an $R(Z)$-module. By \cite[Lemma 4.6]{EdGr:05} 
the localization $G(G,S_\Psi)_{{\mathfrak m}_h} = G(Z,X^h)_{{\mathfrak m}_h}$
at ${\mathfrak m}_h$ is a summand
in the localization $G(G,S_\Psi)_{{\mathfrak m}_\Psi}$, and is independent
of the choice of $h \in \Psi$. This summand is called the {\em central
summand} and is denoted by $G(G,S_\Psi)_{c_\Psi}$ \cite[Definition 4.10]{EdGr:05}.
Because $G$ acts with finite stabilizer, the $G$-equivariant and $Z$-equivariant
Grothendieck groups are supported at a finite number of maximal ideals, so completion
and localization coincide.  Thus, 
$$
G(G,S_\Psi)_{c_\Psi} = G(Z,X^h)_{{\mathfrak m}_h} = \widehat{G(Z,X^h)}_h.
$$

The hypothesis that $G$ acts with finite
stabilizer implies that the projection $f \colon S_\Psi \to
X$ is finite so it
induces a pushforward $f_* \colon G(Z,X^h) \to G(G,X)$.
The nonabelian localization theorem for actions with finite
stabilizer \cite[Theorem 5.1]{EdGr:05} states that
if $\alpha \in G(G,X)_{{\mathfrak m}_\Psi}$ then 
\begin{equation}\label{eqn.oldtheorem}
\alpha = f_*\left(\lambda_{-1}(N_f^*)^{-1} \cap (f^*\alpha)_{c_\Psi}\right)
\end{equation}
where $(f^*\alpha)_{c_\Psi}$ denotes the projection of $f^*\alpha$
to the central summand and $N_f$ is the relative tangent
bundle of the l.c.i. morphism $S_\Psi \to X$.

\begin{remark}
Because the map $f \colon S_\Psi \to X$ is not an embedding,
$N_f^*$ is only a virtual bundle, and the class
$\lambda_{-1}(N_f^*)$ is not defined in $K_0(G,S_\Psi)$.  
Nevertheless, we can define the operation $\cap \lambda_{-1}(N_f^*)^{-1}$ 
on $G(G,S_\Psi)_{c_\Psi}$, as follows.  
Since $f$ factors as the composition of the smooth
projection $G \times_Z X \to X$ with the regular embedding $i \colon
S_\Psi \to G \times_Z X$, we have
$[N_f^*]  = [N_i^*] -  {\mathfrak g^*}/{\mathfrak z^*}$ (here, as usual,
we follow the conventions of Section \ref{ss.background} in identifying
relative tangent bundles).
After localizing $G(G,S_\Psi)$ at ${\mathfrak m}_h$, the endomorphism
$\cap \lambda_{-1}(N_i^*)$ of $G(G,S_\Psi)_{c_\Psi}$ is invertible. Hence
we can define $\cap \lambda_{-1}(N_f^*)^{-1}$ 
as the composition of the inverse of the operation $\cap \lambda_{-1}(N_i^*)$
with multiplication by $\lambda_{-1}({\mathfrak g^*}/{\mathfrak z^*})$.
This point was obscured in \cite{EdGr:05}.
\end{remark}

By abuse of notation
use $i$ to denote both the $Z$-equivariant embedding
$i \colon X^h \to X$ and the $G$-equivariant embedding
$i\colon S_\Psi \to G \times_Z X$; this convention has the advantage
that under our Morita equivalence identifications the
same symbol $i_*$ denotes the proper pushforwards
$G(Z,X^h) \to G(Z,X)$ and $G(G,S_{\Psi}) \to G(G,G \times_Z X)$.
Likewise, the conormal bundles to both of these maps may be denoted
$N_i^*$.
With this convention we can rewrite \eqref{eqn.oldtheorem} as
\begin{equation} \label{eqn.oldtheorem2}
\alpha = f_*\left(\lambda_{-1}({\mathfrak g^*}/{\mathfrak z^*})
(\lambda_{-1}(N_i^*)^{-1} \cap (f^*\alpha)_{c_\Psi }) \right).
\end{equation}
Under our identifications, if $\alpha \in \widehat{G(G,X)}$, then
$f^*(\alpha)_{c_\Psi} = i^!(\alpha)$.  The following result implies that
when $G$ acts on $X$ with finite stabilizers then the nonabelian localization theorem of \cite{EdGr:05} 
is in fact equivalent to Theorem \ref{thm.nonabelian}.

\begin{thm} \label{thm.mapsagree}
Assume that $G$ acts on the smooth space $X$ with finite stabilizer.
If $\beta \in G(G,S_\Psi)_{c_\Psi}$, then
\begin{equation} \label{eqn.mapsagree}
f_*\left(\lambda_{-1}(\mathfrak{g^*/z^*}) \cap \beta\right) = i_!\beta
\end{equation}
(as elements of $G(G,X)_{{\mathfrak m_\Psi}}$).
\end{thm}

\begin{proof}
Since the map $\alpha \mapsto f^*(\alpha)_{c_\Psi} = i^!(\alpha)$ is an 
isomorphism
and $G(G,S_\Psi)_{c_\Psi} = \widehat{G(Z,X^h)}_h$,
we can compare formulas \eqref{eqn.nonabelian}
and \eqref{eqn.oldtheorem} to obtain
\begin{equation}
\label{eqn.temp} f_*(\lambda_{-1}(N_f^*)^{-1} \cap \beta) = i_! (\lambda_{-1}((N_i^*)^{-1}
\cap \beta)).
\end{equation}
But as noted in the remark above, $$\lambda_{-1}(N_f^*)^{-1} \cap \beta
= \lambda_{-1}(N_i^*)^{-1} \cap 
(\lambda_{-1} ({\mathfrak g^*}/{\mathfrak z^*}) \cap \beta). $$
Replacing $\beta$ by
$\lambda_{-1}(N_i^*) \cap \beta$ in \eqref{eqn.temp} yields the
formula of equation \eqref{eqn.mapsagree}.
\end{proof}

\subsection{More on the finite stabilizer case}
The proof of Theorem \ref{thm.mapsagree} used the localization and completion
theorems.  However, if
$G$ is connected and reductive with simply connected commutator subgroup and
$Z$ is a Levi factor  
of parabolic subgroup $P$ of $G$ then, by Proposition \ref{prop.simplelevi}, $Z$ also has a simply 
connected commutator. In this case we can prove a stronger
result (Theorem \ref{thm.mapsagree2}) namely,
that the maps agree even before localizing (or completing).
This proof does not make use of the localization or completion theorems.
Using Theorem \ref{thm.mapsagree2} one could replace some
constructions from \cite{EdGr:05}  with constructions from this paper and
obtain somewhat different proofs of the main results in \cite{EdGr:05}.  
However, Theorem
\ref{thm.mapsagree2} is not needed for the main results of this paper.

Assume that we are still in the finite
stabilizer case.  The proper pushforward $f_*$ is defined without
localizing (or completing); that is, $f_*$ maps $G(G,S_{\Psi})$
to $G(G,X)$.  The map $i_!$ can also be defined without localizing
or completing; we identify $G(G,S_{\Psi})$ with $G(Z,X^h)$ and then
define $i_! = \ind \circ i_*$.  In this case, the equality of Theorem \ref{thm.mapsagree}
holds even without localizing:

\begin{thm} \label{thm.mapsagree2}
Suppose that the connected reductive group $G$ acts with finite stabilizer
on a smooth space $X$.  Let $h \in G$ be an element of the semisimple conjugacy
class $\Psi$, and assume that $Z = {\mathcal Z}_G(h)$ is a Levi
factor of a parabolic subgroup $P$ of $G$.  If $\beta \in G(G,S_{\Psi})$, 
then
\begin{equation} \label{eqn.mapsagree2}
f_*\left(\lambda_{-1}(\mathfrak{g^*/z^*}) \cap \beta\right) = i_!\beta
\end{equation}
as elements of $G(G,X)$.
\end{thm}
\begin{proof}
We use the factorization of the map $S_\psi \to G \times_P X$
from Section 4.1 of \cite{EdGr:05}.  Consider the projections
$$
G \times_Z X \stackrel{p}{\rightarrow} G \times_P X \stackrel{q}{\rightarrow} X.
$$
By \cite[Lemma 4.5]{EdGr:05} the composition $j = p \circ i \colon S_\psi \to G \times_P X$
is also a regular embedding and the finite map $f \colon S_\Psi \to X$ is
the composition of $j$ with the smooth proper projection map $q$.
Thus, 
\begin{equation} \label{e.f*}
f_*\left(\lambda_{-1}(\mathfrak{g^*/z^*}) \cap \beta\right) =
q_*j_*\left(\lambda_{-1}(\mathfrak{g^*/z^*}) \cap \beta\right).
\end{equation}
On the other hand, $i_! \beta = \ind \circ i_* \beta$.
By Proposition \ref{p.altinduction} the right hand side of \eqref{e.f*} equals
\begin{equation} \label{e.q*}
q_*
\left( \lambda_{-1}({\mathfrak g^*}/{\mathfrak p^*}) \cap
(p^*)^{-1} i_*\beta\right)
\end{equation} 
(Note that the symbol $(p^*)^{-1}$ was suppressed in the
statement of Proposition \ref{p.altinduction}.)

The projection $p$ 
is a bundle map with fiber isomorphic to the affine space $P/Z$.
The tangent bundle $T_p$ to $p$ is the bundle induced by Morita
equivalence from the $Z$-module $\mathfrak{p}/\mathfrak{z}$. By abuse
of notation we will write $T_p = \mathfrak{p}/\mathfrak{z}$. Letting
the unipotent radical of $P$ act trivially we may also view $\mathfrak{p}/\mathfrak{z}$
as a $P$-module. It follows that $T_p$ is the pullback of the $G$-bundle
on $G \times_P X$ induced by Morita equivalence from $\mathfrak{p}/\mathfrak{z}$ (viewed
as $Z$-module). Continuing our abuse of notation we will also refer to
this bundle on $G \times_P X$ as $\mathfrak{p}/\mathfrak{z}$.

Claim: If $\beta \in G(G,S_\Psi)$ then 
\begin{equation} \label{e.claim}
i_*\beta = \lambda_{-1}({\mathfrak p^*}/{\mathfrak z^*}) \cap p^* j_*\beta
\end{equation}

We prove the claim.  Let $\rho \colon W \to S_\Psi$ be the $P/Z$-bundle obtained
by base change from the bundle $p \colon G \times_Z X \to G \times_P X$ along the morphism
$j \colon S_\psi \to G \times_P X$. Let $k \colon W \to G \times_Z X$ be
the map obtained  by base change from $j$. The map $i \colon S_\Psi \to G \times_Z X$ induces
a section $s \colon S_\Psi \to W$ of the projection $\rho$.  
In particular we have a cartesian diagram
$$\begin{array}{ccc}
W & \stackrel{k} \to &G \times_Z X\\
\downarrow \rho & & \downarrow p\\
S_\Psi & \stackrel{j} \to & G \times_P X
\end{array}$$
where the horizontal maps are regular embeddings.

Since direct image for finite morphisms commutes with flat pullback of
coherent sheaves, the compositions $p^* \circ j_*$ and $k_* \circ \rho^*$
are equal as maps on $K$-theory. Thus the right hand side of \eqref{e.claim} equals
\begin{equation} \label{e.claimcont}
\lambda_{-1}({\mathfrak p^*}/{\mathfrak z^*}) \cap k_* \rho^*
\beta
\end{equation}
Since the bundle
$\lambda_{-1}(\mathfrak{p^*}/\mathfrak{z^*})$ on $W$ is the pullback of
the corresponding bundle on $G \times_Z X$, 
the the projection formula implies that the right hand side of
\eqref{e.claimcont} equals 
\begin{equation} \label{e.claimcont2}
k_*(\lambda_{-1}({\mathfrak p^*}/{\mathfrak z^*}) \cap \rho^*\beta)
\end{equation}

The normal bundle to the section $s \colon S_\Psi \to W$ is the restriction
of $T_\rho$ to $s$.   Now, $T_{\rho}$ is the pullback of $T_p$,
which is the pullback to $G \times^Z X$ of the
$G$-bundle on $G \times_P X$ induced by 
Morita equivalence from $\mathfrak{p}/\mathfrak{z}$ (viewed
as $P$-module).  By abuse of notation we will refer to this bundle on $G \times_P X$ ,
as well as the various pullbacks of this bundle, as $\mathfrak{p}/\mathfrak{z}$.
Thus, by the self intersection formula $s^*s_* \beta =
\lambda_{-1}(\mathfrak{p^*}/\mathfrak{z^*}) \cap \beta$. Substituting into
\eqref{e.claimcont2}  yields
\begin{equation} \label{e.claimcont3}
k_*(\rho^* s^* s_* \beta)
\end{equation}
Since $s$ is a section of the smooth morphism $\rho$, the composition
$\rho^* \circ s^*$ is the identity. Hence
we conclude that 
$\lambda_{-1}(\mathfrak{p^*}/\mathfrak{z^*}) \cap p^*j_*\beta = k_*s_*\beta = \beta$,
proving our claim.

Substituting the expression for $i_*\beta$ for \eqref{e.claim} into the right
hand side of \eqref{e.q*} we obtain
\begin{equation} \label{e.almostthere}
q_*(\lambda_{-1}({\mathfrak g^*}/{\mathfrak z^*}) \cap \left( \lambda_{-1}(\mathfrak{p^*}/
\mathfrak{z^*}) \cap j_*\beta\right)
\end{equation}
Since all of the terms involving Lie algebras correspond to bundles which
are pulled back from $G \times_P X$, the projection formula implies
that \eqref{e.almostthere} can be rewritten as 
\begin{equation} \label{e.keeppushing}
q_*j_*\left(\lambda_{-1}({\mathfrak g^*}/{\mathfrak p^*}) 
\lambda_{-1}({\mathfrak p^*}/{\mathfrak z^*}) \cap \beta\right).
\end{equation}
Since $f = q\circ j$ and $\mathfrak{g^*}/\mathfrak{z^*} = \mathfrak{g^*}/\mathfrak{p^*} +
\mathfrak{p^*}/\mathfrak{z^*}$ in $R(Z)$ we may rewrite \eqref{e.keeppushing}
as
$f_*(\lambda_{-1}(\mathfrak{g^*}/\mathfrak{z^*}) \cap \beta)$,
completing the proof of the theorem.
\end{proof}

\section{The twisted equivariant Riemann-Roch theorem}
As above, let $\Psi = C_G(h)$ be a semi-simple conjugacy class.
In this section we prove a Riemann-Roch isomorphism
generalizing
the main theorem of \cite{EdGr:00}, which dealt with
the case $h = 1$. This theorem
gives a geometric description of the completion of equivariant
$K$-theory at any maximal ideal in $R(G)$.
\subsection{The equivariant Riemann-Roch theorem}
We recall here the equivariant Riemann-Roch theorem 
of \cite{EdGr:00}.
We begin by recalling some basic facts about
equivariant Chow groups.  Most of these facts are contained
in the paper \cite{EdGr:98}, but here we continue to use
some notation from \cite{EdGr:00}. If $X$ is a $G$-space then
as in \cite{EdGr:00} we denote $CH^n_G(X)$ the ``codimension $n$''
equivariant Chow groups. An element of $CH^n_G(X)$
is represented by a codimension $n$ cycle
on a quotient $X \times_G U$; here $U$ is an open set 
in a representation $V$ such that $G$ acts freely on $U$,
and such that the complement of $U$ has codimension more than $n$ in $V$.  
Less precisely,
but more intuitively, we may view an element of $CH^n_G(X)$
as being represented by a $G$-invariant cycle of codimension $n$
on a product $X \times V$ where $V$ is a representation of $G$.
We denote by $A^k_G(X)$ the ``codimension $n$'' operational
Chow ring. An element $x \in A^k_G(X)$ is an operation
on $CH^*_G(X)$ which increases degree by $k$ and satisfies basic
naturality properties with respect to equivariant morphisms.
The direct sum $A^*_G(X) = \oplus_{k=0}^\infty A^k_G(X)$ is a graded,
commutative ring with multiplication defined by composition of operations.
When $X$ is a smooth algebraic space, then a basic result is that the map 
$A^*_G(X) \to  CH^*_G(X)$, $x \mapsto x \cap [X]$ is an isomorphism,
where $[X] \in CH^n_G(X)$ is the fundamental class.
Thus, when $X$ is smooth we may view the product of operations as 
an intersection product.

The infinite direct product $\prod_{k =0}^\infty A^k_G(X)$ is a
completion of the equivariant operational Chow ring; it acts on the
infinite product of equivariant Chow groups $\prod_{n=0}^\infty
CH^n_G(X)$.  If $E$ is an equivariant vector bundle then we may define
the equivariant Chern character $\ch^G(E)$ and equivariant Todd class
$\Td^G(E)$ as the appropriate formal power series in the Chern classes
of $E$, \cite[Definition 3.1]{EdGr:00}. As such they may be viewed as
elements in the infinite product $\prod_{n=0}^\infty A^n_G(X)$.  
The Chern character and Todd class extend to maps $\ch^G$ and $ \Td^G$
from $K_0(G,X)$ to $\prod_{n=0}^\infty A^n_G(X)$.

We paraphrase from \cite{EdGr:00} the Riemann-Roch theorem for
equivariant Chow groups.
Let $X$ be a $G$-space and 
let $\widehat{G_0(G,X)}$ denote the completion of $G_0(G,X)$ at the
augmentation ideal in $R(G)$---in other words, the
completion at the maximal ideal ${\mathfrak m}_1 
\subset R(G)$.

\begin{thm} \label{thm.rr}
For all separated $G$-spaces $X$ there is a
map $\tau_X^G  \colon G_0(G,X) \to \prod_{n=0}^\infty CH^n_G(X)$
which factors through an isomorphism 
$\widehat{G_0(G,X)} \to \prod_{n=0}^\infty CH^n_G(X)$
and satisfies the following properties.

(a) $\tau_X^G$ is covariant for proper morphisms.

(b) If $\epsilon \in K_0(G,X)$ and $\alpha \in G_0(G,X)$ then
$\tau^G_X(\epsilon \alpha) = \ch^G(\epsilon) \tau^G(\alpha).$

(c) If $X$ is a smooth scheme and either $G$ is
connected or $X$
has a $G$-equivariant ample line bundle, then
$\tau^G(\alpha)= \ch(\alpha) \Td^G(X)$,
where $\Td^G(X)$ is the equivariant Todd class of
the virtual bundle $T_X - {\mathfrak g}$. (Here we implicitly
identify $K_0(G,X)$ and $G_0(G,X)$ as $X$ is a smooth scheme.
\end{thm}

\begin{remark} The notation used here differs somewhat
from that of \cite{EdGr:00}: higher $K$-theory 
is not used in that paper, and the groups
$G^G(X)$ and $K^G(X)$ of that paper
correspond to the groups denoted here by $G_0(G,X)$ and $K_0(G,X)$,
respectively.
\end{remark}

\begin{remark}
The need to subtract
${\mathfrak g}$ follows from the fact that if $G$ acts freely
on $X$ and $X \stackrel{\pi} \to X/G$ is the quotient map, then, by definition,
the Todd class map is the pullback of the Todd
class map  on $X/G$ to the corresponding equivariant theories on $X$. However 
as an element of $K_0^G(X)$, $\pi^*(T_{X/G}) = T_X - {\mathfrak g}$.
\end{remark}

\subsection{Twisting by a central element}
We recall the following construction from \cite{EdGr:05}.  Let $Z$ be
an algebraic group, and let $h$ be a semisimple
element of the center of $Z$.  If $V$ is any representation
of $Z$, then $V$ decomposes as a direct sum
of $H$-eigenspaces
$V_\chi$, where $\chi \in \C^*$.
The assignment $[V]
\to \sum_{\chi} \chi [V_\chi]$ defines an automorphism of
$R(Z)$ which we call twisting by $h$.   It is clear from the
definition that this automorphism takes $\m_h$ to the augmentation
ideal $\m_1$.

More generally, let $X$ be a $Z$-space on which $h$ acts trivially.
Then any $Z$-equivariant coherent sheaf ${\mathcal F}$
decomposes into a direct sum of eigensheaves ${\mathcal F}_\chi$ for the action of $h$. 
Define the twist of ${\mathcal F}$ by $h$ 
to be the virtual coherent sheaf ${\mathcal F}(h) = \sum_{\chi} \chi {\mathcal F}_{\chi}$. 
It is easy to see that twisting by a central 
element induces a natural automorphism of equivariant 
Grothendieck groups. The basic properties of this automorphism
are as follows.
\begin{prop} \label{prop.twist}
Let $X$ be a $Z$-space on which the central element $h$
acts trivially. Then
the assignment ${\mathcal F} \mapsto {\mathcal F}(h)$
induces automorphisms of equivariant 
Grothendieck rings $t_h \colon K_0(Z,X) \to K_0(Z,X)$
and of equivariant Grothendieck groups $t_h \colon G_0(Z,X) \to G_0(Z,X)$
with the following properties.

(a) If $\epsilon \in K_0(Z,X)$ and $ \alpha \in G_0(Z,X)$,
then $t_h(\epsilon \alpha) = t_h(\epsilon) t_h(\alpha)$.

(b) If $f \colon Y \to X$ is an morphism of $Z$-spaces and $h$
acts trivially on $Y$,
then $f^* \circ t_h = t_h \circ f^*$ as maps $K_0(Z,X) \to K_0(Z,Y)$.

(c) $t_h$ commutes with flat pullbacks and proper pushforwards
of Grothendieck groups of equivariant coherent sheaves, provided
the element $h$ acts trivially on the spaces involved.

(d) $t_h$ induces isomorphisms of 
completions 
$\widehat{K_0(Z,X)}_h \to \widehat{K_0(Z,X)}_1$
and
$\widehat{G_0(Z,X)}_h \to \widehat{G_0(Z,X)}_1$.
\end{prop}

\subsection{The twisted equivariant Riemann-Roch theorem}
We can now state the twisted equivariant Riemann-Roch theorem.
For each algebraic $G$-space $X$, and each semisimple conjugacy
class $\Psi$ in $G$,
this theorem gives a way to represent elements of the ${\mathfrak m}_\Psi$-adic
completion of $G_0(G,X)$
by $Z$-equivariant cycles
on $X^h$.   

As usual, we fix $\Psi$ and let $\widehat{G_0(G,X)}$ denote the
${\mathfrak m}_\Psi$-adic
completion of $G_0(G,X)$.

\begin{thm} \label{thm.rrnextgeneration}
For all separated $G$-spaces
$X$, there is an isomorphism $\tau_X^\Psi \colon \widehat{G_0(G,X)}
\to \prod_{i =0}^\infty CH^i_Z(X^h)$ with the following
properties.

(a) $\tau_X^\Psi$ is covariant for proper morphisms.

(b) If 
$\epsilon \in K_0(G,X)$ then 
$$\tau_X^\Psi(\epsilon \alpha)
= \ch^Z(t_h(i^!\epsilon)) \cap
\tau_X^\Psi(\alpha).$$

(c) If $X$ is smooth, then $\ch^Z(t_h(N_i^*))$ is invertible
in the infinite product $\prod_{n=0}^\infty CH^n_Z(X^h)$.  If in addition
$X^h$ is a scheme, and either $Z$ is connected or
$X^h$ has a $Z$-equivariant 
ample line bundle, then
$$
\tau_X^\Psi(\alpha) =  \ch^Z(t_h(i^!\alpha)) \ch^Z(t_h(\lambda_{-1}(N_i^*)))^{-1} \Td^Z(X^h).
$$

(d) The map 
$\tau^\Psi$ is uniquely determined for $G$-schemes by properties (a) and (c).
\end{thm}
\begin{remark}
When $X$ and $X^h$ are smooth schemes then both $G_0(G,X) =K_0(G,X)$
and $\prod_{i=0}^\infty CH^i_Z(X^h)$ have products. The  map $\tau_X^\Psi$ is
not a ring isomorphism but standard properties of the
Chern character imply that the map
$K_0(G,X) \to \prod_{i=0}^\infty CH^i_Z(X^h)$, $\epsilon \mapsto \ch(t_h(i^!\epsilon))$
induces a ring isomorphism $\widehat{G_0(G,X)} \to \prod_{i=0}^\infty CH^i_Z(X^h)$.
\end{remark}
\begin{remark}
We conjecture that $\tau^\Psi$ is also uniquely determined by properties 
(a) and (c) for
all algebraic $G$-spaces. However, the proof of uniqueness uses
an equivariant version of Chow's lemma for
schemes (Proposition \ref{prop.equivchowlemma}),
and we do not know if the corresponding result holds for algebraic spaces.
\end{remark}

\begin{proof}
Define $\tau_X^\Psi$ to be the composition $\tau^Z_{X^h} \circ t_h \circ (i_!)^{-1}$.   
The fact that $\tau_X^\Psi$ is an isomorphism, as well as the covariance for
proper morphisms, follows because these properties hold for
$i_!$  (Theorem \ref{thm.nonabelian}), 
$t_h$ (Proposition \ref{prop.twist}), and $\tau^Z_{X^h}$
(Theorem \ref{thm.rr}).

We now prove part (b). Let $G \subset Q$ be an embedding of $G$ into a
connected reductive group such that $Z_Q =
{\mathcal Z}_Q(h)$ is also connected and reductive, and such that
$\Psi_Q \cap G = \Psi$ .  As in the proof
of Theorem \ref{thm.nonabelian}, let $Y = Q \times_G X$.  In the notation
of that proof, the map
of completions $i_! \colon \widehat{G(Z,X)}_h \to \widehat{G(G,X)}$ is
defined as the composition of maps $\ind_h \circ j_* \circ (\cap \;
\lambda_{-1}(N_e^*)^{-1})$.  Here 
$$\ind_h \colon \widehat{G(Z_Q,Y)} \to \widehat{G(Q,Y)}= \widehat{G(G,X)}$$ 
is the induction map corresponding to the inclusion $Z_Q \subset Q$,
$j_* \colon \widehat{G(Z,X^h)}_h = \widehat{G(Z_Q, Y^h)}_h \to
\widehat{G(Z_Q,Y)}_h$ is the map induced by the inclusion $j \colon
Y^h \to Y$, and $N_e$ is the normal bundle to $Z_Q/Z$ in $Q/G$.  Thus,
if $\epsilon \in K_0(G,X)$, then
\begin{eqnarray*}
(i_!)^{-1}(\epsilon \alpha)  & = & \lambda_{-1}(N^*_e)\cap (j_*)^{-1} \left(
\res_h(\epsilon\alpha ) \right)\\
& = & \lambda_{-1}(N_e^*) \cap  j^*\res_h(\epsilon) \cap
(j_*)^{-1}\left(\res_h\alpha\right)\\
& = & j^!\epsilon  \cap \lambda_{-1}(N_e^*) \cap (j_*)^{-1}\left(\res_h\alpha
\right)\\
& = & i^!\epsilon \cap  (i_!)^{-1}(\alpha).
\end{eqnarray*}
The first equality holds since by Proposition \ref{prop.indres}, 
$\ind_h$
and $\res_h$ are inverse isomorphisms 
of $\widehat{G(Q,Y)} $ and $ \widehat{G(Z_Q,Y)}_h$. The
second equality is a $K$-theoretic consequence of
the projection formula for the direct image functor $j_*$. 
The third equality follows from the definition of $j^!$.  Finally,
the fourth equality follows from the fact that $i^!$ and $j^!$
are equal as maps $K_0(G,X)= K_0(Q,Y) \to K_0(Z_Q,Y^h) = K_0(Z,X^h)$.
Applying Theorem \ref{thm.rr}(b), we see that
$$\tau_{X}^{\Psi} (\epsilon \alpha) = \ch(t_h(i^! \epsilon)) \cap
\tau_{X}^{\Psi}(\alpha),$$ proving (b).

We now prove (c).  
Since multiplication by $\lambda_{-1}(N_i^*)$ is an automorphism
of $\widehat{G_0(Z,X^h)}_h$, it follows that multiplication
by $t_h(\lambda_{-1}(N_i^*))$ is an automorphism of
$\widehat{G_0(Z,X^h)}_1$. Since $\tau_{X^h}^Z \colon \widehat{G_0(Z,X^h)}_1
\to \prod_{i=0}^\infty CH^i_Z(X^h)$ is an isomorphism, and 
$$\tau_{X^h}^Z(t_h(\lambda_{-1}(N_i^*)) \cap \alpha)
 = \ch(t_h(\lambda_{-1}(N_i^*))) \tau_{X^h}(\alpha),$$ it follows
that multiplication by $\ch(t_h(\lambda_{-1}(N_i^*)))$ is an isomorphism,
proving the first assertion.  The formula for $\tau_X^\Psi(\alpha)$ now follows from Theorem \ref{thm.nonabelian} and
Theorem \ref{thm.rr}(c).

The proof of Part (d) uses equivariant Chow envelopes and
will be deferred to Section \ref{sec.envelopeproof}.
\end{proof}

By taking $X$ to be a point, we obtain the following corollary relating
completions of representation rings and Chow rings of classifying stacks.
(By definition, the Chow ring of the classifying stack $BG$ is the
$G$-equivariant Chow ring of a point.)

\begin{cor} \label{cor.rrnextgeneration}
Let $G$ be a complex algebraic group, $\Psi$ the conjugacy
class of a semisimple element $h \in G$, and 
$Z = {\mathcal Z}_G(h)$.
The map $R(G) \to \prod_{i=0}^\infty CH^i(BZ)$ which
takes
$[V]$ to $ \ch(t_h([\res V])$ induces an isomorphism
$\widehat{R(G)} \to \prod_{i=0}^\infty CH^i(BZ)$.
\end{cor}

\section{Equivariant envelopes and proofs of uniqueness}
In this section we use envelope arguments to prove the uniqueness
statements Theorem \ref{thm.nonabelian}(c) and Theorem
\ref{thm.rrnextgeneration}(d), as well as the isomorphism of completions 
stated in Corollary
\ref{cor.nonabelian}(b).

\subsection{Equivariant envelopes}
\begin{definition} \label{def.envelopes}
Let $X$ be a $G$-space.  A
proper $G$-equivariant map $\tilde{X} \to X$ is an 
{\em equivariant envelope} if
for every $G$-invariant subspace $W \subset X$ 
whose connected components are integral, there is
a corresponding $G$-invariant subspace 
$\tilde{W} \subset \tilde{X}$, with integral connected
components, such that the map $\tilde{W} \to W$
is an isomorphism
over a dense open subspace of $W$. 
\end{definition}

If $G$ is connected then every component of
an invariant subspace is invariant, so it would suffice
to require the existence of $\tilde{W}$
for every $G$-invariant subvariety $W$.

By canonical resolution of singularities for algebraic
spaces (cf. \cite{BiMi:97}),
every separated algebraic $G$-space over a field of characteristic
0
has a nonsingular equivariant envelope
\cite[Corollary 2]{EdGr:00}.

The following lemma will be used in the proofs of Theorem \ref{thm.nonabelian}(c)
and Corollary \ref{cor.nonabelian}(b).

\begin{lemma}\label{lem.fixedenvelope}
Let $X$ be a separated $G$-space, and
let $\pi \colon \tilde{X} \to X$ be an equivariant envelope. If
$h \in Z(G)$ then the induced map $\tilde{X}^h \to X^h$ is also
an equivariant envelope.
\end{lemma}
\begin{proof}
Since $h$ is central in $G$, the fixed loci $X^h$ and $\tilde{X}^h$
are $G$-invariant, so the map $\pi \colon \tilde{X}^h \to X^h$ is
$G$-equivariant. 
Let $W \subset X^h$ be any $G$-invariant subspace with integral
connected components, and let $\tilde{W}$ be
a subspace of $\tilde{X}$ as in Definition \ref{def.envelopes}.  
To prove that $\tilde{X}^h$ is an equivariant envelope
it suffices to prove that we can take
$\tilde{W}\subset \tilde{X}^h$.  The restricted map
$\pi \colon \tilde{W} \to W$ is a
$G$-equivariant isomorphism over a dense open subspace of $U$ of $W$.
Replace $\tilde{W}$ by the closure of $\pi^{-1}(U)$.
Since $h$ acts trivially on $\pi^{-1}(U)$ and $\tilde{W}^h$ is closed,
$\tilde{W} \subset \tilde{X}^h$ as desired.
\end{proof}

\subsubsection{Equivariant Chow envelopes} 
Let $X$ be an algebraic space with a $G$-action. 
An equivariant envelope $\tilde{X} \to X$ is an {\it equivariant Chow
envelope} if $\tilde{X}$ has a $G$-linearized ample line bundle.
In this section we prove that nonsingular
Chow envelopes exist in the case where $X$ is a separated
scheme (Proposition \ref{prop.chowenvelopes}).  As a first step, we
need an equivariant version of Chow's lemma, which extends
a result of Sumihiro to disconnected groups.  
Once we establish this version of Chow's lemma, Proposition
\ref{prop.chowenvelopes} will follow 
from equivariant resolution of singularities and Noetherian
induction.

\begin{prop}[Equivariant Chow's Lemma] 
\label{prop.equivchowlemma} Let $X$ be a reduced and separated $G$-scheme.
There exists a proper $G$-map $X'
\stackrel{\pi} \to X$ such that $X'$ has a $G$-linearized
ample line bundle, and such that there is a dense open $G$-invariant 
subspace $U \subset X$ such that $\pi^{-1}(U) \to U$ is an isomorphism.
\end{prop}
\begin{proof}[Proof of Proposition \ref{prop.equivchowlemma}]
When $G$ is connected this is proved by Sumihiro in \cite{Sum:74}.
In general, let $G_0$ be the identity component of $G$.  
Let $X_1 \stackrel{\pi} \to X$ be a proper
$G_0$-map with $X'$ normal and quasi-projective and $\pi$ an isomorphism over a
dense $G_0$-invariant open set $U$ in $X$.  
By replacing $U$ by the intersection of all $\sigma U$, where
$\sigma$ runs over a set of coset representatives of $G_0$ in $G$,
we may assume that $U$ is $G$-invariant.  Then the inverse image
$U_1$ of $U$ admits a $G$-action (since it is isomorphic to $U$).
To produce a $G$-equivariant Chow cover we
apply a construction used in the proof of \cite[Theorem 3]{Sum:74}.

If $Y$ is any $G_0$-space, let $IY$ denote the set of all functions
$f: G \to Y$ satisfying $f(g g_0) = g_0^{-1}f(g)$ for all $g \in G$ and
$g_0 \in G_0$.  As a space, $IY$ is isomorphic to $Y^r$ by the
map taking $f$ to $ (f(\sigma_1) \ldots, f(\sigma_r))$, where
$\sigma_1, \ldots, \sigma_r$ is a collection of coset representatives
of $G_0$ in $G$.  (This isomorphism depends on the choice
of coset representatives.)  
The group $G$ acts on $IY$ by the rule
$(g \cdot f)(h) = f(g^{-1} h)$.   

If $Y$ has a $G_0$-equivariant ample line bundle, then
$IY$ has a $G$-equivariant ample line bundle.  This follows because
a $G_0$-equivariant embedding of $Y$ into  $\Pro^N$
induces a $G$-equivariant embedding of $IY$ into $I \Pro^N \simeq (\Pro^N)^r$.
By Kambyashi's theorem \cite{Kam:66}, 
$(\Pro^N)^r$ has a $G$-linearized ample line bundle and this restricts to 
a $G$-equivariant ample line bundle on $IY$.

The projection $P \colon IY \to Y$,
which takes $f$ to $f(e)$ is $G_0$-equivariant. 
If in addition $Y$ is a $G$-space then the map $\Delta \colon Y \to IY$
which takes $y \in Y$ to the function $f \colon G \to Y; g \mapsto g^{-1}y$
is a $G$-invariant section of $P$. In this case we denote the image of the section by
$\Delta_Y$.

We now return to the proof the proposition.

The proper $G_0$-equivariant map $\pi: X_1 \to X$ induces a proper
$G$-equivariant map $I \pi : IX_1 \to IX$.  Let $X'$ denote the inverse
image of $\Delta_X$ under this map.  Then the map $X' \to \Delta_X \simeq X$ is proper
and $G$-equivariant.  Since $U_1$ is a $G$-space, $\Delta_{U_1} \simeq U_1$
embeds
in $X'$ and maps isomorphically to $\Delta_U \simeq U$.  It  follows that $X'
\to X$ is an isomorphism over the dense open set $U$.   Finally, there is a
$G$-equivariant ample line bundle on $IX_1$, and this restricts to a $G$-equivariant
ample bundle on $X'$.
\end{proof}

To prove Proposition
\ref{prop.chowenvelopes}, we also need the following lemma,
which is stated without proof in \cite{PoVi:94}.

\begin{lemma}\label{lem.linearize}
Let $G$ be an algebraic group acting on a 
normal quasi-projective variety $Y$. Then 
there exists a $G$-equivariant embedding 
$Y \subset \Pro^N$
such that the $G$-action is the restriction of the action of 
$\PGL_{N+1}$ on $\Pro^N$; i.e., $Y$ has a $G$-linearized ample line bundle.
\end{lemma}
\begin{proof}[Proof of Lemma \ref{lem.linearize}]
If $G$ is connected then our statement is \cite[Theorem 1]{Sum:74}.
Otherwise let $G_0$ be the identity component of $G$, 
and let $Y \hookrightarrow \Pro^N$ be an $G_0$-equivariant projective
embedding. As in the proof of Proposition \ref{prop.equivchowlemma}, 
since $Y$ is a $G$-space we have $Y \simeq \Delta_Y \subset IY$.
The composite embedding of $Y$ into $I \Pro^N \simeq (\Pro^N)^r$
 is $G$-equivariant. Since $(\Pro^N)^r$ is
projective, it has a $G$-linearized ample line bundle by Kambayashi's
theorem \cite{Kam:66}. Pulling this ample bundle back to $X$ gives us
our desired $G$-linearized ample line bundle.
\end{proof}

\begin{prop} \label{prop.chowenvelopes}
If $X$ is a separated $G$-scheme, then there exists
an equivariant Chow envelope $\tilde{X} \to X$ with $\tilde{X}$
non-singular.
\end{prop}
\begin{proof}
The map $X_{red} \to X$ is an equivariant envelope,
as is the normalization map. Thus we may assume that $X$ is normal.
We will show that
there is a non-singular quasi-projective variety
$X'$ and a proper $G$-equivariant morphism 
$X' \stackrel{\pi}\to X$ such that $\pi$ is an isomorphism
over a dense $G$-equivariant open subspace $U$ of $X$.
This implies the proposition, since by Noetherian induction
we may assume there is an equivariant Chow envelope $\tilde{Y}$
of $Y = X-U$, so $\tilde{X} = X' \cup \tilde{Y} \to X$ is the
desired envelope.

By Proposition \ref{prop.equivchowlemma} there is a quasi-projective
scheme $X_1$ and a proper $G$-equivariant morphism $X_1 \to X$
which is an isomorphism over a dense open subspace of $X$.
By canonical resolution of singularities there is a canonical (hence
$G$-equivariant) resolution of singularities $X' \to X_1$. Since
$X'$ is obtained by a sequence of blowups of the quasi-projective
scheme $X_1$, $X'$ is
quasi-projective. The composite map $X' \to X$ is
an isomorphism over a dense open subspace of $X$. Finally, since
$X'$ is non-singular it has a
a $G$-linearized ample line bundle by Lemma \ref{lem.linearize}.
\end{proof}

\subsubsection{Envelopes and completions of equivariant $K$-theory}
A crucial property of (equivariant) envelopes
is the following:
\begin{lemma} \label{lem.envpushforward}
Let $X$ be a $G$-space and $\pi \colon \tilde{X} \to X$
a birational envelope. Then $\pi_* \colon CH^*_G(\tilde{X}) \to CH^*_G(X)$
is surjective.
\end{lemma}
\begin{proof}
This follows from \cite[Lemma 3]{EdGr:98} and the corresponding
fact for non-equivariant envelopes \cite[Lemma 18.3(6)]{Ful:84}
\end{proof}

Using Riemann-Roch and this lemma, we can prove an 
analogous result on completions 
of equivariant Grothendieck groups:

\begin{prop} \label{prop.kenvelope}
Let $\tilde{X}\stackrel{\pi} \to X$ be a $G$-equivariant
envelope. Then for any $h \in Z(G)$, the proper pushforward
$\pi_* \colon \widehat{G_0(G,\tilde{X}^h)}_h \to \widehat{G_0(G,X^h)}_h$ is 
surjective.
\end{prop}
\begin{proof}
By Theorem \ref{thm.rr} and Proposition \ref{prop.twist}
the diagram
$$\begin{array}{ccc}
\widehat{G_0(G,\tilde{X}^h)}_h & \stackrel{\tau^G \circ t_h}\to &
\prod_{n=0}^\infty CH^n_G(\tilde{X}^h)\\
\pi_* \downarrow & & \pi_* \downarrow\\
\widehat{G_0(G,X^h)}_h & \stackrel{\tau^G \circ t_h} \to &
\prod_{n=0}^\infty CH^n_G(X^h)
\end{array}$$
commutes and the horizontal arrows are isomorphisms. By Lemma 
\ref{lem.envpushforward}, $\pi_*$ is surjective as a map of 
equivariant Chow groups. Hence $\pi_*$ must also be surjective
as a map of completed equivariant 
Grothendieck groups.
\end{proof}

\subsection{Proofs of Theorem \ref{thm.nonabelian}(c), 
Corollary \ref{cor.nonabelian}(b) and Theorem \ref{thm.rrnextgeneration}(d)} 
\label{sec.envelopeproof}
Suppose that there is an assignment 
$j_X \colon \widehat{G_0(Z,X^h)}_h \to \widehat{G_0(G,X)}$ 
for every $G$-space $X$ which is covariant
for proper morphisms and satisfies equation \eqref{eqn.nonabelian}
for smooth $X$. Equation \eqref{eqn.nonabelian}
implies that 
$$j_X \alpha = (i^!)^{-1}\left( \lambda_{-1}(N_i^*)^{-1} \cap \alpha\right)
= i_! \alpha$$ when $X$ is smooth.
We wish to show that $j_X = i_!$ for arbitrary $X$. Let 
$\pi \colon \tilde{X} \to X$ be a $G$-equivariant envelope
with $\tilde{X}$ non-singular. By Proposition \ref{prop.kenvelope},
$\pi_* \colon \widehat{G_0(Z,\tilde{X}^h)}_h \to \widehat{G_0(Z,X^h)}_h$
is surjective. By covariance, $j_X(\pi_*\alpha) = \pi_*j_{\tilde{X}} \alpha$.
Since $j_{\tilde X} = i_!$ as maps $\widehat{G_0(Z,\tilde{X}^h)}_h \to
\widehat{G_0(G,\tilde{X})}$, it follows from covariance of $i_!$ that
$\pi_*j_{\tilde{X}} \alpha = i_!\pi_* \alpha$. Hence $j_X  = i_!$
completing the proof of Theorem \ref{thm.nonabelian}(c).

To prove part (b) of Corollary \ref{cor.nonabelian} we use a
similar argument. First observe that since $i_!$ is an isomorphism
which is covariant for proper morphisms, Proposition
\ref{prop.kenvelope} implies that the pushforward
$\pi_*\colon  \widehat{G_0(G,\tilde{X})} \to \widehat{G_0(G,X)}$
is surjective. Consider, for each $G$-space $X$, the map
$(i_*)^{-1} \circ \res_h \colon \widehat{G_0(G,X)} \to
\widehat{G_0(Z,X^h)}_h$. This map is covariant for proper morphisms.
When $X$ is smooth, equation \eqref{eqn.nonabelian} implies
that $(i_!)^{-1}$ is the map
$$
\alpha \mapsto \lambda_{-1}(N_i^*) \cap i^* ( \res_h \alpha) = (i_*)^{-1} \circ \res_h (\alpha),
$$
i.e., that $(i_*)^{-1} \circ \res_h = (i_!)^{-1}$. Using the same
argument as in the previous paragraph, we conclude that
$(i_*)^{-1} \circ \res_h = (i_!)^{-1}$ as maps
$\widehat{G_0(G,X)} \to \widehat{G_0(Z,X^h)}_h$
for arbitrary $X$.
Since $i_*$ and $i_!$ are isomorphisms, it follows that
$\res_h \colon \widehat{G_0(G,X)} \to \widehat{G_0(Z,X)}_h$
is also an isomorphism. 

We conclude with a proof of Theorem 
\ref{thm.rrnextgeneration}(d),  the uniqueness
of the functor $\tau^\Psi$.
If $X$ is a $G$-scheme, then by
Proposition \ref{prop.equivchowlemma}
there is a non-singular equivariant Chow envelope
$\tilde{X} \to X$. If there is a functor $\tau'$ such that
$\tau'$ satisfies properties (a) and (c), then $\tau'= \tau^\Psi$
on the smooth quasi-projective scheme $\tilde{X}$. Arguing as above
shows that $\tau' = \tau^{\Psi}$ for arbitrary $X$.
%\bibliographystyle{amsmath}
%\bibliography{refs}
%\bibliography{jabbrev,refs}
\def\cprime{$'$}

\end{document}